\documentclass[journal]{IEEEtran}

\ifCLASSINFOpdf
\else
\fi
\usepackage{algorithm}
\usepackage[noend]{algorithmic}
\usepackage{placeins}
\usepackage{multicol}
\usepackage{pgfplots}
\usepackage{subcaption}
\usepackage{amsmath}
\DeclareMathOperator{\atantwo}{atan2}
\usepackage{url}

\usepackage{color}

\pgfplotsset{compat=1.17}

\begin{document}
%



\title{The Promise of EV-Aware Multi-Period OPF Problem: Cost and Emission Benefits}

%
%
%
\author{Sezen Ece Kayac{\i}k, Burak Kocuk, Tu\u{g}\c{c}e Y\"{u}ksel
\thanks{
S. E. Kayac{\i}k is with the Faculty of Economics and Business, University of Groningen, Groningen, Netherlands (e-mail: s.e.kayacik@rug.nl).

B. Kocuk and T. Y\"{u}ksel are with the Faculty of Engineering and Natural Sciences, Sabanc{\i} University, Istanbul, Turkey (e-mail: burak.kocuk, tugce.yuksel@sabanciuniv.edu).
}
}

\newcommand{\cI}{\mathcal{I}}
\newcommand{\cS}{\mathcal{S}}
\newcommand{\cT}{\mathcal{T}}
\newcommand{\cB}{\mathcal{B}}
\newcommand{\cL}{\mathcal{L}}
\newcommand{\cG}{\mathcal{G}}

\maketitle

\begin{abstract}
 In this paper, we study the Multi-Period Optimal Power Flow problem (MOPF) with electric vehicles (EV) under emission considerations. We integrate three different real-world datasets: household electricity consumption, marginal emission factors, and EV driving profiles. We present a systematic solution approach based on second-order cone programming to find globally optimal solutions for the resulting nonconvex optimization problem. To the best of our knowledge, our paper is the first to propose such a comprehensive model integrating multiple real datasets and a promising solution method for the EV-aware MOPF problem. Our computational experiments on various instances with up to 2000 buses demonstrate that our solution approach leads to high-quality feasible solutions with provably small optimality gaps. In addition, we show the importance of coordinated EV charging to achieve significant emission savings and reductions in cost. In turn, our findings can provide insights to decision-makers on how to incentivize EV drivers depending on the trade-off between cost and emission.
\end{abstract}


%
\IEEEpeerreviewmaketitle

\section{Introduction}
\IEEEPARstart{M}{any} countries, including the United States, grant  tax credits and other incentives to increase electric vehicle (EV) usage aiming to reduce emissions \cite{AFDC}. 
Since EV charging demand is expected to comprise a major part of the electricity consumption in the future, charging strategies will become more critical from different perspectives such as network security, financial parties, environmental policy, and EV drivers.

EVs create a considerable opportunity to reduce greenhouse gas emissions from transportation compared to conventional vehicles. However, the additional electricity generation to charge EVs may cause greenhouse gas emissions depending on the source of electricity generation \cite{yuksel2016effect}. Therefore, it is also important to consider the additional amount of emission generated by the power plants to correctly judge the environmental benefits of EVs.

Integrating EV charging stations to power systems changes the dynamics of conventional power distribution and poses considerable challenges to network operations under an uncoordinated charging scenario. These challenges include increased network losses, overloaded transmission lines, unbalanced voltage, and increased peak load, resulting in a significant increase in electricity generation cost and the risk of damage to power system security. Recently, various coordinated charging strategies of EVs have been developed to overcome economic and security issues by flattening the load profiles \cite{clement2009impact}, filling the load valley \cite{gan2012optimal} and minimizing the stress on transformers \cite{masoum2012distribution}.  


 
With the increased integration of electric vehicles, Optimal Power Flow (OPF) studies that consider large-scale integration of EVs have become increasingly important to protect network security. The classical OPF problem considers single-period planning and optimizes operations of the power grid with respect to a certain objective and network constraints. The interconnection of the energy storage of EVs with the power network makes it necessary to solve a multi-period version of the OPF problem (in short, MOPF). 


In recent years, several studies have focused on the OPF problem that considers the integration of EVs with limited scope. Some papers study single-period OPF and analyze the effects of EV penetration by resolving the problem with the inclusion of  EVs \cite{judd2008evaluation}.
Some studies do not consider the effects of EV charging on power system security and ignore operational limits such as limits of line flow, generation and, voltage \cite{acha2010effects}.
The authors in \cite{tang2016model} do not consider the physical limitations of batteries and assume that EVs can be fully charged during a single period. Many existing studies incorporate only a small number of EVs \cite{6847105} or charging stations \cite{fan2017admm}. The work in \cite{shi2018model} considers only economic aspects without environmental impacts. The authors of \cite{7384482} consider both operational cost and emission
objectives; however, they do not guarantee to produce globally optimal solutions and carry out tests only on small-size networks. Study by \cite{fan2017admm} do not consider the behavior of EV users and randomly generate arrival times. 

Systematic treatments to solving large-scale MOPF problems with EVs are limited in the literature. Since Alternating Current (AC) OPF is a nonlinear, nonconvex, and NP-hard problem, a globally optimal solution cannot be guaranteed in general. The integration of EVs requires a multi-period setting and additional constraints, increasing the complexity and computational effort even more. 
Primarily, simulation-oriented \cite{6485952,yang2014improved} and heuristic approaches \cite{valentine2011intelligent,alonso2014optimal,yang2014self} are developed, however, these methods fail to provide a quality measure for global optimality. A comprehensive overview of heuristic approaches can be found in \cite{yang2015computational}.

 Convex relaxations of the single-period OPF problem have been widely studied  \cite{zohrizadeh2020survey}, however very few studies focus on convex relaxations of the MOPF problem. The authors in \cite{li2012optimization} model a demand response problem in radial distribution networks as an AC OPF problem ignoring reactive power. This optimization problem is convexified and extended to an AC MOPF problem. Reference \cite{gopalakrishnan2013global} solves the AC MOPF problem using a Semidefinite Programming relaxation, and reference \cite{jabr2014minimum} proposes a Second-Order Cone Programming (SOCP) relaxation of the AC MOPF problem for distribution networks with photovoltaics.  The study by \cite{huang2016sufficient} formulates the AC MOPF with the planning of EV charging as an SOCP problem.

In this paper, we study an EV-aware AC MOPF problem with an integrated perspective to preserve system security, capture the trade-off between cost and emission, and satisfy the charging needs of EV owners simultaneously. We develop a convex programming framework to find globally optimal solutions for large-scale networks in a reasonable time. We consider both Grid to Vehicle (G2V) and Vehicle to Grid (V2G) concepts, where the latter allows EVs to supply power back to the grid. To implement V2G, we integrate an \textit{aggregator} concept in order to control a group of EVs to function as a source of electricity. We analyze the impacts of this bidirectional flow on the grid considering several other aspects, including penetration level, driving profiles, and battery characteristics. The main contributions of this paper  can be summarized as follows:

\begin{itemize}
\item 

We introduce a new EV-aware MOPF formulation for the joint problem of OPF and EV charging by bringing together test cases and real data for electricity consumption, marginal emission factors, and EV driving profiles.
\item 

We develop a convex optimization framework to solve the resulting challenging problem. Specifically, we introduce an SOCP model in which an emission constraint is imposed to restrict emissions caused due to EV charging from power plants. EV charging constraints and their effects on the power system are also considered.
\item 

We demonstrate the effectiveness of the proposed approach on three test cases with up to 2000 buses from PGLIB-OPF. Our algorithm can solve large-scale test cases in reasonable computation time and offer guarantees of global optimality with very small optimality gaps.
\item 

Our extensive computational experiments suggest that through coordinated charging of EVs, marginal emission can be decreased while keeping the cost virtually constant, and the integration of the V2G concept leads to cost savings, despite assuming hourly electricity prices are constant.
\end{itemize}

The remainder of this paper is organized as follows. Section \ref{sec:model} describes the problem,  provides an alternative formulation and its SOCP relaxation. Section \ref{sec:MOPFSolApproach} introduces our approach to solve this  optimization problem. Section \ref{sec:MOPFInput} introduces the input data used in the formulation. Section \ref{sec:MOPFCompExperiments} presents the experimental setup and the computational results. The concluding remarks are made in Section \ref{sec:MOPFConclusions}. 

\section{Mathematical Model}\label{sec:model}
\subsection{Formulation}
\label{sec:MINLP}
In this section, we present our mathematical programming formulation and its SOCP relaxation. We formulate an MOPF problem that coordinates the charging of EVs. The model minimizes the total cost of generation over a finite planning horizon while determining the optimal schedule of power generation and EV charging under emission limits. We take into account both network constraints and EV charging constraints.
 
 We consider power station emissions associated with EV charging only. Therefore, we utilize marginal emissions, which occur to satisfy the additional demand due to EVs. The EV charging constraints are formulated based on actual driving profiles that provide arrival and departure times of each trip and energy requirements for each period. We assume that only the real power is consumed and fed back to the grid by EVs. All EVs connected to the same bus are aggregated as an entity. In the formulation, we will use the term EV to refer to a group of EVs charged at the same bus.
Note that our formulation considers aggregate charging load and requires perfect knowledge of the user's driving behavior and electricity consumption. 

Consider a power network $\mathcal{N} = (\mathcal{B}, \mathcal{L})$, where $\mathcal{B}$ denotes the node set, i.e., the set of buses, and $\mathcal{L}$ denotes the edge set, i.e., the set of transmission lines. Let   $\mathcal{G} \subseteq \mathcal{B}$ denote the set of generators, $\mathcal{D} \subseteq \mathcal{B}$ be the set of buses with real power load and $\mathcal{T} = \{1,...,T\}$ represent the set of time periods.

Rest of the parameters are given as follows:
\begin{itemize}
\item
For each bus $i \in \mathcal{B}$ and time $t \in \mathcal{T}$;
$p_{it}^d$ and $q_{it}^d$ are the real and reactive power load, $\underline V_{i}$ and $\overline V_i$ are the bounds on the voltage magnitude,  $\delta(i)$ is the set of its neighbors, $b_{ii}$ is the shunt susceptance, $\overline{a}_{it},\overline{b}_{it}$ are the maximum allowable charging and discharging rates of EV, $\underline{s}_{it}, \overline{s}_{it}$ are the minimum energy requirement and maximum capacity of EV battery, $c_{it}$ is the energy requirement of EV, $\tilde p_{it}^g$ is the reference power generation without EVs, $\eta_i$ the charging efficiency of EV and $I_i$ is the initial battery state of charge of EV for bus $i$.
\item
For each generator located at bus $i \in \mathcal{G}$ and time $t \in \mathcal{T}$; active and reactive outputs must be in the intervals $ [\underline p_i, \overline p_i]$ and $[\underline q_i , \overline q_i]$.

\item
For each line $(i,j) \in \mathcal{L}$ and time $t \in \mathcal{T}$; $G_{ij}, \textcolor{black}{G_{ji}} $  and $B_{ij}, \textcolor{black}{B_{ji}}$  are the elements of the conductance and susceptance matrices, $\overline S_{ij}$ is the apparent power flow limit and $\overline\theta_{ij}$ is the bound on the phase angle.
\item
For each time $t \in \mathcal{T}$; $e_{t}$ is the marginal emission parameter and $\overline{E}$ is the upper bound on the total emission over the planning horizon.
\end{itemize}
\ \  We define the following decision variables:
\begin{itemize}
\item 
For each bus $i \in \mathcal{B}$ and time $t \in \mathcal{T}$, $|V_{it}|$ and $\theta_{it}$ are the voltage magnitude and phase angle, $a_{it},b_{it}$ are the charging and discharging power of EV, $s_{it}$ is the stock variable of EV.

\item
For each generator located at bus $i \in \mathcal{G}$ and time $t \in \mathcal{T}$, $p_{it}^g$ and $q_{it}^g$ are the real and reactive power output.

\item
For each line $(i,j) \in \mathcal{L}$ and time $t \in \mathcal{T}$, $p_{ijt}$ and $q_{ijt}$ are the real and reactive power flow.
\end{itemize}
We present the following formulation of the EV-aware MOPF problem: 
\begin{subequations}\label{MMNLP}
\allowdisplaybreaks\begin{align}
\min  &\hspace{0.25em}   \sum_{i\in \mathcal{G}} \sum_{t\in\mathcal{T}} f(p_{it}^g) \label{Mobj} \\
  \mathrm{s.t. \ }   & \text{For each bus } i \in \mathcal{B} \text{ and time t} \in \mathcal{T},\notag \\ 
  &\hspace{0.25em} p_{it}^g - p_{it}^d -a_{it} +\eta_i b_{it} =  \ g_{ii} |V_{it}|^2 +  \sum_{j\in\delta(i)} p_{ijt}   \label{MactiveAtBus} \\
  & \hspace{0.25em} q_{it}^g - q_{it}^d = -b_{ii} |V_{it}|^2 + \sum_{j\in\delta(i)} q_{ijt} \label{MreactiveAtBus} \\
  & \text{For each line } (i,j) \in \mathcal{L} \text{ and time t} \in \mathcal{T}, \notag\\ 
  & \hspace{0.25em}  p_{ijt} = G_{ij}|V_{it}|^2 + |V_{it}|  |V_{jt}|[G_{ij} \cos(\theta_{it} - \theta_{jt}) \label{MactiveFlowAtLine} \\
  &\hspace{4.3cm}-B_{ij} \sin(\theta_{it} - \theta_{jt}) ]   \nonumber \\
  & \hspace{0.25em} q_{ijt} = - B_{ij}|V_{it}|^2   - |V_{it}| |V_{jt}|[B_{ij}\cos(\theta_{it} - \theta_{jt}) \label{MreactiveFlowAtLine} \\
  & \hspace{4.3cm} +G_{ij} \sin(\theta_{it} - \theta_{jt}) ]   \nonumber \\
  & \hspace{0.25em} \underline V_i \le |V_{it}| \le \overline  V_i   \hspace{3cm} i \in \mathcal{B},t \in \mathcal{T}    \label{MvoltageAtBus} \\
  & \hspace{0.25em}  \underline{p}_i  \le p_{it}^g \le \overline{p}_i  \hspace{3.5cm} i \in \mathcal{G},t \in \mathcal{T}  \label{MactiveAtGenerator} \\
  & \hspace{0.25em} \underline{q}_i  \le q_{it}^g \le \overline{q}_i  \hspace{3.5cm} i \in \mathcal{G},t \in \mathcal{T}    \label{MreactiveAtGenerator}\\
   & \hspace{0.25em} p_{ijt}^2 + q_{ijt}^2  \le \overline{S}_{ij}^{2}   \hspace{2.4cm} (i,j) \in \mathcal{L},t \in \mathcal{T} \label{MpowerOnArc}\\
   & \hspace{0.25em} |\theta_{it} - \theta_{jt}| \le \overline\theta_{ij}  \hspace{2.5cm} (i,j) \in \mathcal{L},t \in \mathcal{T}  \label{MphaseAngle} \\
   & \hspace{0.25em} s_{it}+\eta_i a_{it}-b_{it}-c_{it}=s_{i(t+1)}  \hspace{0.7cm} i \in \mathcal{B},t \in \mathcal{T}  \ \label{Mstock constraint} \\
   & \underline s_{it} \leq s_{it} \leq \overline s_{it} \hspace{3.3cm} i \in \mathcal{B},t \in \mathcal{T} \label{Mstock limit} \\
   & s_{i0}=s_{iT}=I_i \hspace{4.3cm}  i \in \mathcal{B} \label{Mstock initial} \\
   & 0 \leq a_{it} \leq \overline a_{it} \hspace{3.45cm} i \in \mathcal{B},t \in \mathcal{T}  \label{Mcharging limita} \\
   & 0 \leq b_{it} \leq \overline b_{it} \hspace{3.5cm} i \in \mathcal{B},t \in \mathcal{T}  \label{Mcharging limitb} \\
   & \sum_{i \in \mathcal{B}}\sum_{t \in \mathcal{T}} e_{t} (p^g_{it} - \tilde p_{it}^g ) \leq \overline E. \label{Memission}
\end{align}
\end{subequations}
Here, the objective function \eqref{Mobj} aims to minimize the total cost of power generation. Constraint \eqref{MactiveAtBus} ensures real power flow balance at bus $i$, while considering EV charging and discharging power. Constraint \eqref{MreactiveAtBus} ensures reactive power flow balance at bus $i$. Constraints \eqref{MactiveFlowAtLine} and \eqref{MreactiveFlowAtLine} represent the real and reactive power flow, respectively. Constraint \eqref{MvoltageAtBus} enforce bus voltage magnitude to maintain a level under acceptable limitations. Constraints \eqref{MactiveAtGenerator} and \eqref{MreactiveAtGenerator} limit real and reactive power outputs of generator $i$ (we set $\underline{p}_i = \overline{p}_i =\underline{q}_i = \overline{q}_i = 0$ for $ i \in \mathcal{B} \setminus  \mathcal{G}$). Constraint \eqref{MpowerOnArc} satisfy transmission capacity limitations of line $(i,j)$. Constraint \eqref{MphaseAngle} sets restrictions on phase angle. Constraint \eqref{Mstock constraint} ensures balance between power supply and EV demand. Constraint \eqref{Mstock limit} controls the load of EV battery in a specified range between the minimum charging requirement and the maximum battery capacity in order to satisfy the charging demand of the EV. Constraint \eqref{Mstock initial}  set the battery state of charge at the beginning and end of the planning horizon. Constraints \eqref{Mcharging limita} and \eqref{Mcharging limitb} limit the charging and discharging rate of EV. If an EV is connected to the grid at time $t$, its charging/discharging rate should be between zero and the maximum limit. Otherwise, maximum allowable charging and discharging limits $(\overline{a}_{it},\overline{b}_{it})$ are set to zero. Constraint \eqref{Memission} sets an upper limit on the total amount of marginal emission allowed, according to the difference between generation with and without EVs.

\subsection{Alternative Formulation}\label{sec:MOPFAlternative}

We now present an alternative formulation for the mathematical model \eqref{MMNLP} motivated by \cite{kocuk2016strong}. Let us first define the new decision variables: 
\begin{subequations}
\allowdisplaybreaks\begin{align}
  & c_{iit} := |V_{it}|^2  & i & \in \mathcal{B}, t \in \mathcal{T} \notag  \\
  & c_{ijt} := \ |V_{it}||V_{jt}|\cos(\theta_{ij} - \theta_{jt}) & (&i,j) \in \mathcal{L}, t \in \mathcal{T} \notag \\
   & s_{ijt} := -|V_{it}||V_{jt}| \sin(\theta_{it} - \theta_{jt}) & (&i,j) \in \mathcal{L}, t \in \mathcal{T}. \notag 
\end{align}
\end{subequations}

To eliminate the nonlinearities, we rewrite the constraints \eqref{MactiveAtBus}--\eqref{MvoltageAtBus} for each $t \in \mathcal{T}$, replacing the newly defined variables  as follows:
\begin{subequations}\label{MRNLP}
\allowdisplaybreaks
\begin{align}
   & p_{it}^g - p_{it}^d - a_{it}+ \eta_i b_{it}=  \ g_{ii} c_{iit} +  \sum_{j\in\delta(i)} p_{ijt} \hspace{0.3cm} i \in \mathcal{B}  \ \ \  \label{MRactiveAtBus} \\
  &  q_{it}^g - q_{it}^d=   -  {b_{ii}} {c_{iit}}  + \sum_{j\in\delta(i)} q_{ijt} \hspace{2.3cm} i \in \mathcal{B} \label{MRreactiveAtBus} \\
  &  p_{ijt} = G_{ij} c_{iit}   +  G_{ij}c_{ijt}    - B_{ij}s_{ijt} \hspace{1.5cm} (i,j) \in \mathcal{L}
    \label{MRactiveFlowAtLine}  \\
  & q_{ijt} = \  -B_{ij} c_{iit}  -  B_{ij}c_{ijt}  - G_{ij} s_{ijt} \hspace{1.2cm} (i,j) \in \mathcal{L}
  \label{MRreactiveFlowAtLine}  \\
  & \underline V_i^2 \le c_{iit} \le \overline  V_i^2 \hspace{4.6cm} i \in \mathcal{B}. \label{MRvoltageAtBus} 
\end{align}
\end{subequations}

We define the following consistency constraints which preserves the trigonometric relation between the variables $c_{iit}, c_{ijt}$ and $s_{ijt}$:
\begin{subequations}\label{MCNLP}
\allowdisplaybreaks\begin{align}
  & c_{ijt}^2+s_{ijt}^2  = c_{iit} c_{jjt}  & (&i,j) \in \mathcal{L}, t \in \mathcal{T} 
  \label{MRconsistency1}  \\
  & \theta_{jt}- \theta_{it} = \atantwo(s_{ijt},c_{ijt}) & i& \in \mathcal{B}, t \in \mathcal{T}. \label{MRconsistency2} 
\end{align}
\end{subequations}

The alternative formulation minimizes the objective function \eqref{Mobj} subject to  constraints \eqref{MactiveAtGenerator}--\eqref{Memission}, \eqref{MRNLP} and \eqref{MCNLP}.

\subsection{SOCP Relaxations}\label{sec:MOPFSOCP}
The feasible region of the alternative NLP formulation is nonconvex due to constraints \eqref{MCNLP}. We eliminate the constraint \eqref{MRconsistency2} and relax the constraint \eqref{MRconsistency1} as follows: 
\begin{equation}\label{MRconsistency relaxed}
    c_{ijt}^2+s_{ijt}^2  \le c_{iit} c_{jjt}    \ \ \  (i,j) \in \mathcal{L}, t \in \mathcal{T}.
\end{equation}
Then, an SOCP relaxation of the EV-aware MOPF problem is obtained as minimizing
\eqref{Mobj} subject to \eqref{MactiveAtGenerator}--\eqref{Memission}, \eqref{MRNLP} and \eqref{MRconsistency relaxed}.  

\noindent
\section{Solution Approach}\label{sec:MOPFSolApproach}

Due to its nonconvexity, obtaining a globally optimal solution of the original problem \eqref{MMNLP} is quite challenging. In this section, we propose a solution approach based on its SOCP relaxation. The exactness of this convex relaxation would guarantee that the resulting solution is also globally optimal for the original problem. Although the SOCP relaxation is rarely exact in practice, we still make use of it in two aspects: Firstly, it provides a lower bound for the optimal value of the original problem \eqref{MMNLP}. Secondly, we utilize the optimal solution of the relaxation to guide the local solver to obtain a feasible solution, hence, an upper bound, for the original nonconvex problem. This outlined procedure provides both lower and upper bounds for the optimal value of problem \eqref{MMNLP}, from which we can compute a quality measure for global optimality.

We introduce Algorithm \ref{Algorithm1} in order to explore the trade-off between cost and marginal emission. We aim to approximate the Pareto frontier following a multi-objective optimization approach. Since we are interested in marginal emissions related to the additional generation to satisfy EV demand, we first find generation values without EVs, denoted by $\tilde p_{it}$. To find a lower bound on marginal emission (LBE), we solve the problem under the emission minimization objective. To find an upper bound on marginal emission (UBE), we first solve the problem under the cost minimization objective and then calculate the marginal emission. We create a list of emission bounds varying between LBE and UBE. For each fixed bound $\overline{E}$, we repeatedly perform the optimization algorithm to obtain non-dominated solutions, which  make up the  Pareto frontier.

We start the optimization algorithm by solving the SOCP relaxation from Section \ref{sec:MOPFSOCP}. Since the original problem is challenging to solve, we fix the charging and discharging power variables $a_{it}$ and $ b_{it}$ to their optimal values obtained from the SOCP relaxation in the MOPF model. Consequently, the EV charging constraints \eqref{Mstock constraint}--\eqref{Memission} are no longer needed. These simplifications enable us to decompose the remaining problem into $T$ subproblems with respect to time periods. We solve these subproblems, denoted by $SP_t$, in parallel and sum up their objective function values to obtain an upper bound (UB). We also calculate corresponding marginal emission values. Since we eliminate emission constraint \eqref{Memission} from the original problem, marginal emission values of original and SOCP problems may not match. Therefore, we solve SOCP relaxation one more time fixing $\overline{E}$ to marginal emission from the original problem to obtain a lower bound (LB). The optimality gap is computed as $\%\text{Gap}=100\times(1-\text{LB}/\text{UB})$.


\begin{algorithm}
\caption{}
\label{Algorithm1}
\begin{algorithmic}[1]
\STATE Solve $\min $\{\eqref{Mobj} : \eqref{MactiveAtBus}--\eqref{MphaseAngle}\} where $(a_{it},b_{it})$ equal to zero to obtain $ p_{it}^g{^*}$ 
\STATE Set $ \tilde p_{it}^g = p_{it}^g{^*} $ 
\STATE 
LBE=$\min \{\sum_{i \in \mathcal{B}}\sum_{t \in \mathcal{T}} e_{t} (p^g_{it} - \tilde p_{it}^g )$ :  \eqref{MactiveAtGenerator}--\eqref{Mcharging limitb}, \eqref{MRNLP}, \eqref{MRconsistency relaxed}\}
\STATE Solve $\min$ \{ \eqref{Mobj} : \eqref{MactiveAtGenerator}--\eqref{Mcharging limitb}, \eqref{MRNLP}, \eqref{MRconsistency relaxed} \} to obtain $(p^g_{it}{^*})$ 
\STATE Compute UBE=$\sum_{i \in \mathcal{B}}\sum_{t \in \mathcal{T}} e_{t} (p^g_{it}{^*} - \tilde p_{it}^g )$
\STATE Generate a list of size  $n$, which consists of numbers spaced 
between LBE and UBE as $\{\rho^1,\rho^2,...,\rho^n \}$
\FOR{$k=1 \ to \ n$}
\STATE Set $\overline{E} = \rho^k$
\STATE Solve $\min $\{ \eqref{Mobj} :\eqref{MactiveAtGenerator}--\eqref{Memission}, \eqref{MRNLP}, \eqref{MRconsistency relaxed} \} to obtain  $(a_{it}^*,b_{it}^*)$
\STATE Fix $(a_{it}^*,b_{it}^*)$ and decompose $\min $\{ \eqref{Mobj} : \eqref{MactiveAtBus}--\eqref{MphaseAngle} \} into $\mathcal{T}$ sub-problems
\FORALL{$t \in \mathcal{T}\ (\text{in parallel})$} 
\STATE Solve $SP_t$ to obtain UB$_{kt}$
\STATE Calculate \text{Emission}$_{kt}= \sum_{i \in \mathcal{B}}\sum_{t \in \mathcal{T}} e_{t} (p^g_{it} - \tilde p_{it}^g )$
\ENDFOR
\STATE UB$_k = \sum_{t \in \mathcal{T}} \text{UB}_{kt}$
\STATE Emission$_k = \sum_{t \in \mathcal{T}} \text{Emission}_{kt}$
\STATE Set $\overline{E} = \text{Emission}_k$
\STATE LB$_k = \min $\{ \eqref{Mobj} :\eqref{MactiveAtGenerator}--\eqref{Memission}, \eqref{MRNLP}, \eqref{MRconsistency relaxed} \}
\STATE Compute  $\%\text{Gap}=100\times(1-\text{LB}_k/\text{UB}_k)$ 
\ENDFOR
\STATE Plot coordinates  $\{ (\text{UB}_{k}, \text{Emission}_k) : k=1 \ ,\dots, \ n\}$
    \end{algorithmic}
\end{algorithm}

\section{Input} \label{sec:MOPFInput}

In this section, we provide information on input datasets utilized for constructing the optimization model. We aim to use realistic datasets in order to suggest practical solutions to cope with possible grid challenges. Most of the researchers use synthetic grid data in OPF studies since access to confidential information on power grids is restricted, and actual grid data is not publicly available. In our study, we use a  synthetic but realistically created test case and integrate real datasets to increase the applicability of the optimization model. In particular, our model draws on three real datasets: hourly electricity consumption, hourly marginal emission factors, and EV driving profiles. We select the regions where a realistic OPF test case and real datasets are available. In the following sections, we explain these datasets and describe their variations across different regions. In Section \ref{Experimental Setup}, we will discuss in detail how to integrate these datasets into our model formulation.


\subsection{OPF Test Case}
The first important piece of information is the OPF test case. We test our algorithm on 200-bus, 500-bus, and 2000-bus Texas A\&M University (TAMU) instances from the Power Grid Library (PGLIB-OPF) \cite{birchfield2016grid}, geographically situated in Illinois (IL), South Carolina (SC), and Texas (TX), respectively. These synthetic test cases do not match any actual grid; they are constructed based on the statistical characteristics of the actual grid. The creators of these test cases initially situate substations in a specified region and determine loads and generators of these substations in such a way that generation and load profiles are similar to real profiles. Then, they link these substations with an automated line placement process based on realistic choices. They also consider additional complexities such as voltage control and transient stability to make the test cases more realistic. For detailed information on how to generate these test cases, see \cite{birchfield2016grid}. The main reason we prefer to use TAMU instances is their similarity in size, complexity, and characteristics to real networks and the availability of real datasets given below for the same regions.

In the test cases used in this study, substations are situated in the central part of the IL, northwestern part of the SC, and the whole state of TX. 
%

\subsection{Electricity Demand}
We retrieve the hourly electricity demand data available from Energy Information Administration \cite{eia.gov2020}. We match each test case with the demand data in the corresponding region. The hourly demand data for Illinois is not available; instead, we consider the Midcontinent Independent System Operator's demand. For SC and TX, we use demand data for the SC Public Service Authority and regional demand of Texas, respectively. 

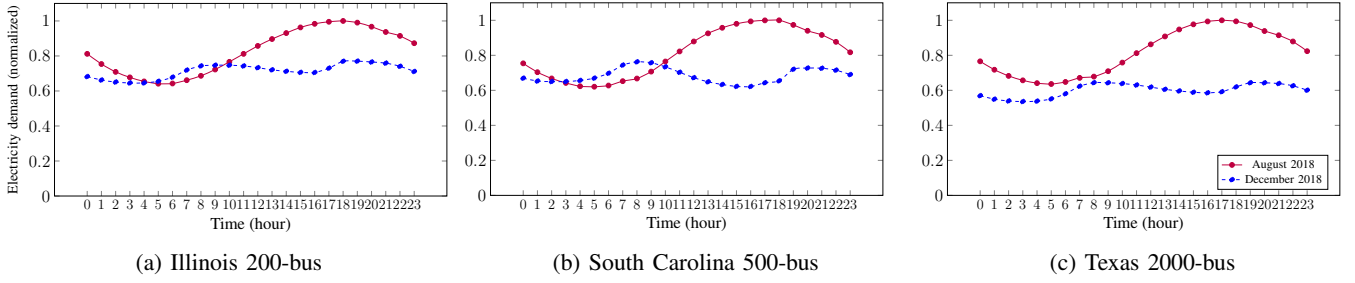
\begin{figure*}[ht]
\allowdisplaybreaks
\begin{subfigure}{0.33\textwidth}
\begin{tikzpicture}[scale=0.45, transform shape]
        \begin{axis}[
            xlabel = \Large Time (hour),
            ylabel = 
                \large Electricity demand  (normalized),
            title style={font=\large},
            xtick = {0,1,2,3,4,5,6,7,8,9,10,11,12,13,14,15,16,17,18,19,20,21,22,23},
            x tick label style={font=\large},
            ytick = {0, 0.2, 0.4, 0.6, 0.8, 1,1.2},
            y tick label style={font=\Large},
            ymin = 0,
            x = 0.42cm,
            legend pos = south east,
            xmajorgrids=false,
            grid style=dashed]
        \addplot[smooth, thick, color = purple, mark = *]
            plot coordinates {
                (0.0,0.8121085414022059)
(1.0,0.7538775716411893)
(2.0,0.7091147781791285)
(3.0,0.6776742976072324)
(4.0,0.6545525266223263)
(5.0,0.6415661455206558)
(6.0,0.6428082062804971)
(7.0,0.6615230933320028)
(8.0,0.6865543663001871)
(9.0,0.7225254930548608)
(10.0,0.7664440386702635)
(11.0,0.8120289788195693)
(12.0,0.8573840706238962)
(13.0,0.8966526261825664)
(14.0,0.930745192842349)
(15.0,0.962990138643542)
(16.0,0.9836277894348086)
(17.0,0.99531906805488)
(18.0,1.0)
(19.0,0.9904082890908243)
(20.0,0.9672290566826952)
(21.0,0.9369643347184432)
(22.0,0.9146956509831263)
(23.0,0.8728324720858358)
                
            };
            \addplot[smooth, thick, color = blue, mark = *,style=dashed]
            plot coordinates {
                (0.0,0.682222625247958)
(1.0,0.6624822640251073)
(2.0,0.6509722104036794)
(3.0,0.6457520218380512)
(4.0,0.6459332468366947)
(5.0,0.6553525734762008)
(6.0,0.6789649790613228)
(7.0,0.7191175624319265)
(8.0,0.7429774978199752)
(9.0,0.7472164145330854)
(10.0,0.7468937445010735)
(11.0,0.7428846748068992)
(12.0,0.7329349309498159)
(13.0,0.7206911339527535)
(14.0,0.7123238023454713)
(15.0,0.7062416844752368)
(16.0,0.7055344619604815)
(17.0,0.7303491470141265)
(18.0,0.7695911817119176)
(19.0,0.7703558656498962)
(20.0,0.765467187845254)
(21.0,0.7589828373603711)
(22.0,0.7408779287830363)
(23.0,0.7112673876117821)
                
            };
        \end{axis}
    \end{tikzpicture}
    \caption{Illinois 200-bus}\label{fig:IllinoisConsumption}
\end{subfigure}
\begin{subfigure}{0.33\textwidth}

\begin{tikzpicture}[scale=0.45, transform shape]
        \begin{axis}[
            xlabel = \Large Time (hour),
            title style={font=\large},
            xtick = {0,1,2,3,4,5,6,7,8,9,10,11,12,13,14,15,16,17,18,19,20,21,22,23},
            x tick label style={font=\large},
            ytick = {0, 0.2, 0.4, 0.6, 0.8, 1,1.2},
            y tick label style={font=\Large},
            ymin = 0,
            x = 0.42cm,
            legend pos = south east,
            xmajorgrids=false,
            grid style=dashed]
        \addplot[smooth, thick, color = purple, mark = *]
            plot coordinates {
                (0.0,0.7541647658449874)
(1.0,0.7035819345907597)
(2.0,0.6679592254086112)
(3.0,0.6422863345341283)
(4.0,0.6233934779530903)
(5.0,0.6207034875489625)
(6.0,0.6274599254353538)
(7.0,0.6523148075319732)
(8.0,0.6674479699224466)
(9.0,0.7070191445515897)
(10.0,0.7647123597980149)
(11.0,0.8222246692570279)
(12.0,0.8791549340087149)
(13.0,0.9254432191791597)
(14.0,0.9574242162060124)
(15.0,0.9800924979156508)
(16.0,0.9938334722899527)
(17.0,0.9996067265491042)
(18.0,1.0)
(19.0,0.9737843917632808)
(20.0,0.9402696282779344)
(21.0,0.916350736994447)
(22.0,0.8770233919048593)
(23.0,0.8172851547137756)
                
            };
            \addplot[smooth, thick, color = blue, mark = *,style=dashed]
            plot coordinates {
                (0.0,0.6693750098318364)
(1.0,0.652999103336532)
(2.0,0.6499473013575799)
(3.0,0.6503484402774937)
(4.0,0.6564441787663798)
(5.0,0.6696424357784455)
(6.0,0.697549119854017)
(7.0,0.7441205619091067)
(8.0,0.7626752033223742)
(9.0,0.7573974736113516)
(10.0,0.7340763579732259)
(11.0,0.7031257373877204)
(12.0,0.6722380405543583)
(13.0,0.6486101716245339)
(14.0,0.6333904890748635)
(15.0,0.622378832449779)
(16.0,0.6218518460255785)
(17.0,0.6432695181613679)
(18.0,0.6542182510343092)
(19.0,0.7201938051566015)
(20.0,0.7270839560162973)
(21.0,0.7260535795749501)
(22.0,0.7153250798345105)
(23.0,0.6898330947474398)
            };
        \end{axis}
    \end{tikzpicture}
    \caption{South Carolina 500-bus} \label{fig:SouthCarolinaConsumption}
\end{subfigure}
\begin{subfigure}{0.33\textwidth}
    \begin{tikzpicture}[scale=0.45, transform shape]
        \begin{axis}[
            xlabel = \Large Time (hour),
            title style={font=\large},
            xtick = {0,1,2,3,4,5,6,7,8,9,10,11,12,13,14,15,16,17,18,19,20,21,22,23},
            x tick label style={font=\large},
            ytick = {0, 0.2, 0.4, 0.6, 0.8, 1,1.2},
            y tick label style={font=\Large},
            ymin = 0,
            x = 0.42cm,
            legend pos = south east,
            xmajorgrids=false,
            grid style=dashed]
        \addplot[smooth, thick, color = purple, mark = *]
            plot coordinates {
                (0.0,0.7656734755770473)
(1.0,0.717730068989997)
(2.0,0.6824580863270027)
(3.0,0.6570668969289155)
(4.0,0.6410689546990574)
(5.0,0.6361488305809974)
(6.0,0.6478732391437999)
(7.0,0.671178079471927)
(8.0,0.6785099279973025)
(9.0,0.7092105498838523)
(10.0,0.7585490145152778)
(11.0,0.8124146601561827)
(12.0,0.8627032618759453)
(13.0,0.9078742255284746)
(14.0,0.9477134103035115)
(15.0,0.976952443311779)
(16.0,0.9936042532655187)
(17.0,1.0)
(18.0,0.9941136265029153)
(19.0,0.9725869106947028)
(20.0,0.9387572860108281)
(21.0,0.914534013927428)
(22.0,0.8790359337092802)
(23.0,0.8236005895086889)
                
            };
            \addplot[smooth, thick, color = blue, mark = *,style=dashed]
            plot coordinates {
                (0.0,0.5695223583075716)
(1.0,0.5494307232056017)
(2.0,0.539403358460965)
(3.0,0.535895962422901)
(4.0,0.5383924085073599)
(5.0,0.5511178903907263)
(6.0,0.5803959075642389)
(7.0,0.6234514161739084)
(8.0,0.6437374434979677)
(9.0,0.6430981287605361)
(10.0,0.638754427962091)
(11.0,0.6303248292534182)
(12.0,0.6182573737583352)
(13.0,0.6056493589411028)
(14.0,0.5963751365839411)
(15.0,0.5890214577992389)
(16.0,0.5857936980273275)
(17.0,0.592024158512619)
(18.0,0.619454398605226)
(19.0,0.6432727708351514)
(20.0,0.6429661087079521)
(21.0,0.6388090036104084)
(22.0,0.6259499408109163)
(23.0,0.6006980487382899)
                
            };
            
        \addlegendentry{August 2018} 
        \addlegendentry{December 2018}  
        \end{axis}
    \end{tikzpicture}
    \caption{Texas 2000-bus}  \label{fig:TexasConsumption}
\end{subfigure}
\caption{Hourly electricity demand. }  \label{fig:Demand}
\end{figure*}

Seasonal variations in electricity demand may result in differences in the operation scheduling of the network. Therefore, the numerical experiments are conducted using data both for a summer's day and a winter's day to demonstrate the influences of the variations on the model results. To generate an \textit{average} day, we first calculate the hourly averages of electricity consumption for a certain month. Then, these hourly averages are normalized by the maximum consumption 
and plotted in Fig. \ref{fig:Demand} for August 2018 and December 2018. We will explain how to process the normalized data to determine the real and reactive power demand of the grid in more detail in Section \ref{Experimental Setup}.


\begin{figure*}[ht]
\allowdisplaybreaks
\begin{subfigure}{0.33\textwidth}
\begin{tikzpicture}[scale=0.45, transform shape]
        \begin{axis}[
            xlabel = \Large Time (hour),
            ylabel = \Large Emission factor (kg/MWh),
            title style={font=\large},
            xtick = {0,1,2,3,4,5,6,7,8,9,10,11,12,13,14,15,16,17,18,19,20,21,22,23},
            x tick label style={font=\large},
            ytick = {400, 600,800,1000},
            y tick label style={font=\large},
            ymax = 1000,
            ymin = 400,
            x = 0.42cm,
            legend pos = north east,
            xmajorgrids=false,
            grid style=dashed]
        \addplot[smooth, thick, color = purple, mark = *]
            plot coordinates {
                (0.0,770.27428553102)
(1.0,783.810372195843)
(2.0,852.962561711628)
(3.0,804.707329093425)
(4.0,778.765680165325)
(5.0,765.157973019088)
(6.0,722.443968724287)
(7.0,750.668374846066)
(8.0,732.602749104127)
(9.0,634.119702915538)
(10.0,599.980771249603)
(11.0,572.362628283221)
(12.0,567.243482687214)
(13.0,599.788018517273)
(14.0,592.902917929997)
(15.0,616.009461864139)
(16.0,653.565730060798)
(17.0,619.823775497369)
(18.0,588.073224778202)
(19.0,579.570011488906)
(20.0,574.338362730963)
(21.0,588.709162151409)
(22.0,650.083113469772)
(23.0,716.977099237854)
            };
            \addplot[smooth, thick, color = blue, mark = *,style=dashed]
            plot coordinates {
                (0.0,709.004340854909)
(1.0,750.699857843924)
(2.0,751.926814496804)
(3.0,792.04197934194)
(4.0,798.274217214963)
(5.0,712.469689547748)
(6.0,654.856562268859)
(7.0,663.758816477344)
(8.0,709.914031747643)
(9.0,703.629116471986)
(10.0,712.813449107387)
(11.0,694.132800498388)
(12.0,692.47935144006)
(13.0,753.665304787468)
(14.0,757.660724173347)
(15.0,795.23401409578)
(16.0,692.624296006482)
(17.0,634.22091099509)
(18.0,670.096278302777)
(19.0,719.476505032267)
(20.0,738.546819940168)
(21.0,633.433786592572)
(22.0,710.736570650539)
(23.0,692.476755409863)
            };
\end{axis}
\end{tikzpicture}
\caption{Illinois 200-bus} \label{fig:IllinoisEmission}
\end{subfigure}
\begin{subfigure}{0.33\textwidth}
 
\begin{tikzpicture}[scale=0.45, transform shape]
        \begin{axis}[
            xlabel = \Large Time (hour),
            title style={font=\large},
            xtick = {0,1,2,3,4,5,6,7,8,9,10,11,12,13,14,15,16,17,18,19,20,21,22,23},
            x tick label style={font=\large},
            ytick = {400, 600,800,1000},
            x tick label style={font=\large},
            y tick label style={font=\large},
            ymax = 1000,
            ymin = 400,
            x = 0.42cm,
            legend pos = north east,
            xmajorgrids=false,
            grid style=dashed]
        \addplot[smooth, thick, color = purple, mark = *]
            plot coordinates {
                (0, 633.594353627174)
                (1, 557.8867247)
                (2, 544.1783852)
                (3, 482.1955296)
                (4, 507.6517323)
    (5, 649.7493863)
    (6, 673.1335105)
    (7, 654.5255434)
    (8, 736.0004681)
    (9, 699.3363716)
    (10, 743.4179282)
    (11, 695.9633562)
    (12, 606.8589041)
    (13, 622.802878)
    (14, 624.1942136)
    (15, 646.7370177)
    (16, 615.0357622)
    (17, 645.1252728)
    (18, 611.9039997)
    (19, 674.6535286)
    (20, 643.3024136)
    (21, 722.0774547)
    (22, 655.102518)
    (23, 628.2557644)
            };
            \addplot[smooth, thick, color = blue, mark = *,style=dashed]
            plot coordinates {
                ( 0 , 662.5303621 )
    (1, 691.2929018)
    (2, 728.0680072)
    (3, 634.969026)
    (4, 642.943483)
    (5, 638.2600538)
    (6, 605.4654505)
    (7, 600.1026594)
    (8, 621.7378489)
    (9, 619.2837739)
    (10, 632.4281045)
    (11, 656.0101095)
    (12, 668.5226922)
    (13, 636.5127004)
    (14, 694.826565)
    (15, 627.2384492)
    (16, 691.9821867)
    (17, 715.4160674)
    (18, 646.2681375)
    (19, 636.7213932)
    (20, 658.7953432)
    (21, 654.164608)
    (22, 647.112126)
    (23, 654.2292527)
            };
        \end{axis}
    \end{tikzpicture}
    \caption{South Carolina 500-bus}  \label{fig:SouthCarolinaEmission}
\end{subfigure}
\begin{subfigure}{0.33\textwidth}
\begin{tikzpicture}[scale=0.45, transform shape]
        \begin{axis}[
            xlabel = \Large Time (hour),
            title style={font=\large},
            xtick = {0,1,2,3,4,5,6,7,8,9,10,11,12,13,14,15,16,17,18,19,20,21,22,23},
            x tick label style={font=\large},
            ytick = {400, 600,800,1000},
            y tick label style={font=\large},
            ymax = 1000,
            ymin = 400,
            x = 0.42cm,
            legend pos = north east,
            xmajorgrids=false,
            grid style=dashed]
        \addplot[smooth, thick, color = purple, mark = *]
            plot coordinates {
(0.0,658.0774566207)
(1.0,714.79009594349)
(2.0,680.12345979579)
(3.0,672.6200492032)
(4.0,667.90888590459)
(5.0,657.10410506715)
(6.0,596.81187797621)
(7.0,670.81342137243)
(8.0,635.71652467779)
(9.0,643.94971277013)
(10.0,646.32714029392)
(11.0,586.1674527884)
(12.0,546.20349939342)
(13.0,529.94551190652)
(14.0,530.77042723863)
(15.0,531.15605998905)
(16.0,530.23811507914)
(17.0,533.63640443644)
(18.0,516.75708734496)
(19.0,538.52089953792)
(20.0,536.10889825315)
(21.0,516.14037391982)
(22.0,532.42537291376)
(23.0,626.0441406735)
            };
            \addplot[smooth, thick, color = blue, mark = *,style=dashed]
            plot coordinates {
                (0.0,594.18539326521)
(1.0,649.23435840918)
(2.0,667.82814526675)
(3.0,684.66465412424)
(4.0,653.83107194181)
(5.0,569.21187685798)
(6.0,533.66102642276)
(7.0,560.57236895625)
(8.0,574.49414311071)
(9.0,580.73052205383)
(10.0,598.84176242298)
(11.0,580.26340821945)
(12.0,533.34280582541)
(13.0,546.48087272473)
(14.0,585.50835415492)
(15.0,595.07847768758)
(16.0,618.98193809689)
(17.0,582.72167875015)
(18.0,563.28551422809)
(19.0,634.78302611049)
(20.0,653.63306542639)
(21.0,611.23274878968)
(22.0,574.19028261344)
(23.0,558.67418280659)
            };
        \addlegendentry{Summer 2018} 
        \addlegendentry{Winter 2018} 
        \end{axis}
    \end{tikzpicture}
     \caption{Texas 2000-bus}  \label{fig:TexasEmission}
\end{subfigure}
\caption{Marginal $\text{CO}_2$ emission factors.}\label{fig:Emission}
\end{figure*}
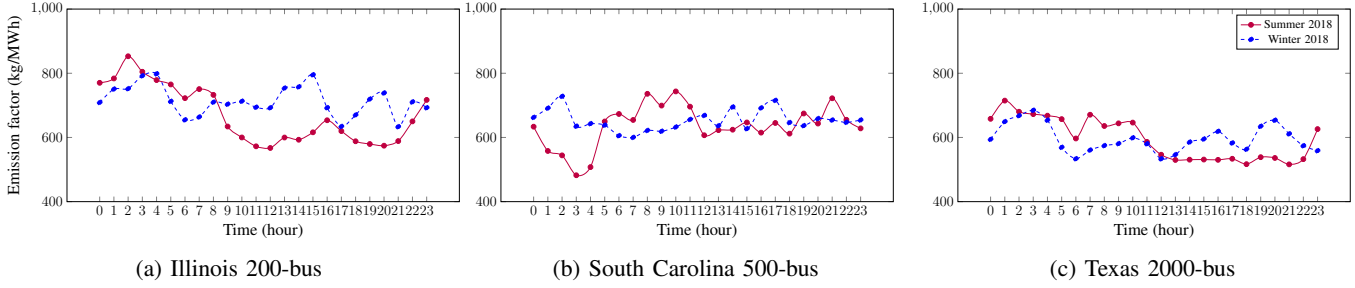

Hourly electricity consumption follows a similar trend across all regions. The total consumption in summer is higher than in winter. In winter, hourly electricity consumption varies slightly and peaks in the morning and evening, while it peaks in the afternoon in summer.

\subsection{Emission Factors} \label{sec:marginalemission}
 We use the marginal emissions that occur from power stations to satisfy the additional demand for EV charging. We retrieve marginal emission factors estimates, available from the Climate and Energy Decision Making Center  \cite{cedm2020}. We only take into account $\text{CO}_2$ emissions. For the winter and summer, the hourly marginal $\text{CO}_2$ emission factors in kg/MWh are shown in Fig. \ref{fig:Emission}. 

The marginal emission factors, as shown in Fig. \ref{fig:Emission}, vary widely across regions, seasons, and time of the day.

\subsection{Driving Profiles}\label{sec:Drivingprofiles}

 Our model also considers driving profiles to determine charging profiles. We use the National Household Travel Survey (NHTS) \cite{NHTS} dataset, a household survey that tracks trends in personal travel and includes all transportation modes. The NHTS dataset contains household, person, vehicle, and travel day trip files. We use vehicle and travel day trip files. The vehicle file includes records for each household vehicle, while the travel day trip file includes records for each personal trip. Regional data is available for IL, SC, and TX. We select trips made by passenger cars and eliminate nonhousehold vehicles. 
  We electrify the gasoline vehicles in the NHTS dataset to investigate how their electricity consumption will affect the OPF problem.
  Then we obtain driving profiles over 24-hour that include energy requirement at time $t$, which is explained in detail in Algorithm~\ref{DrivingProfileAlgorithm}.

\begin{algorithm}[ht]
\caption{Driving Profile}
\label{DrivingProfileAlgorithm}
\textbf{Input:} Trip index $r \in [0,R] $, vehicle index $v \in [0,V]$\\
From vehicle and travel day trip file import $D_r$: Total distance travelled during trip $r$, $S_r$: Start time of trip $r$, $E_r$: End time of trip $r$, $\mathcal{S}_v$: Set of trips for vehicle $v$, $c^{avg}$: Average energy consumption per mile, $\overline{B}$: Battery capacity \\
\textbf{Output:} Driving Profile 
\begin{algorithmic}[1]
\STATE Select passenger cars (Car, SUV, Van, Pickup Truck)
\STATE Delete nonhousehold vehicles
\IF{$D_r \ge \overline{B} / c^{avg} $}
\STATE Delete trip $r$
\ENDIF
\STATE Create a drive duration matrix $\Delta$ of size $T \times R$
\FORALL{$t \in [0,T], r \in [0,R]$}
\IF{$S_r \le t -1$ \textbf{and} $E_r \ge t$}
\STATE $\Delta_{rt} = 60$
\ENDIF
\IF{$S_r \ge t $ \textbf{or} $E_r \le t-1$}
\STATE $\Delta_{rt} = 0$
\ELSE
\STATE $\Delta_{rt} = \min(t,E_r)-\max(t-1,S_r)$
\ENDIF
\ENDFOR
\STATE Calculate average speed $S_r^{avg}= D_r / (E_r -S_r)$
\STATE Create an EV energy consumption matrix $\Omega$ of size $T \times V$
\FORALL{$t \in [0,T], v \in [0,V]$}
\STATE $\Omega_{vt} = \sum_{r \in \mathcal{S}_v} S_r^{avg} \Delta_{tr} c^{avg}  $
\ENDFOR
\STATE Driving profile of vehicle $v$ $\Omega_v=\{ \Omega_{vt} : t \in [0,T] \}$
    \end{algorithmic}
\end{algorithm}

To obtain the energy requirement of vehicles at time $t$, we make some assumptions. Nissan Leaf has the highest sales in the moderate cost electric car segment in US \cite{INSIDEEVs}. Therefore, we assume all of the EVs are similar to Nissan Leaf. The standard battery of Nissan Leaf has a usable capacity of 32 kWh, that is 80\% of its 40 kWh capacity \cite{nissanleaf}; we set the battery capacity $\overline{B}$ to 32. We assume the energy consumption is constant and set the corresponding parameter $c^{avg}$ to average energy consumption value of 0.3 kWh/mile \cite{fueleco}. Lastly, we also eliminate trips that require more than a full battery. After filtering the data, we calculate energy requirements at each time $t$ for vehicle $v$ as outlined in Algorithm \ref{DrivingProfileAlgorithm}. 

\begin{figure}[ht]
    \centering
\begin{center}
    \begin{tikzpicture}[scale=0.56]
        \begin{axis}[
            xlabel = \Large Time (hour),
            ylabel = \large Percentage of energy requirement (\%),
            xtick = {0,1,2,3,4,5,6,7,8,9,10,11,12,13,14,15,16,17,18,19,20,21,22,23},
            x tick label style={font=\large},
            ytick = {0, 2, 4, 6, 8, 10,12},
            y tick label style={font=\Large},
            legend pos=north west,
            x = 0.5cm,
            legend pos = north west,
            xmajorgrids=false,
            grid style=dashed]
            
        \addplot[smooth, thick, color = red,style=dashed]
            plot coordinates {
                (0.0,0.0)
(1.0,0.0)
(2.0,0.0)
(3.0,0.3679)
(4.0,0.2393)
(5.0,3.7638)
(6.0,5.4854)
(7.0,8.8682)
(8.0,5.9555)
(9.0,2.9876)
(10.0,6.303699999999999)
(11.0,7.2285)
(12.0,7.2269)
(13.0,5.0355)
(14.0,4.4415000000000004)
(15.0,7.1729)
(16.0,9.7406)
(17.0,8.791599999999999)
(18.0,7.4106000000000005)
(19.0,2.7474)
(20.0,1.1923)
(21.0,2.1563)
(22.0,2.7824999999999998)
(23.0,0.10189999999999999)
                
            };
\addplot[smooth, very thick, color = blue, style=densely dotted]
            plot coordinates {
                (0.0,0.0)
(0.0,0.0)
(1.0,0.1408)
(2.0,0.2053)
(3.0,0.0)
(4.0,0.7898000000000001)
(5.0,1.9743)
(6.0,5.3558)
(7.0,9.3103)
(8.0,5.21)
(9.0,4.9896)
(10.0,4.9131)
(11.0,4.4068000000000005)
(12.0,6.361999999999999)
(13.0,4.6028)
(14.0,6.8671)
(15.0,7.4675)
(16.0,8.1958)
(17.0,10.1655)
(18.0,6.6551)
(19.0,3.9815000000000005)
(20.0,2.1316)
(21.0,2.8681)
(22.0,2.4572)
(23.0,0.95)
                
            };
            
\addplot[smooth, thick, color = orange]
            plot coordinates {
                (0.0,0.0)
(1.0,0.0834)
(2.0,0.1497)
(3.0,0.0)
(4.0,1.0517)
(5.0,3.8434000000000004)
(6.0,6.6511000000000005)
(7.0,7.4666)
(8.0,6.2)
(9.0,4.8803)
(10.0,4.2523)
(11.0,4.611400000000001)
(12.0,5.382)
(13.0,4.4805)
(14.0,5.0926)
(15.0,7.2146)
(16.0,8.2621)
(17.0,9.476700000000001)
(18.0,7.0665000000000004)
(19.0,5.0657)
(20.0,3.1969)
(21.0,2.6662000000000003)
(22.0,2.0669)
(23.0,0.8392000000000001)
                
            };
            \addlegendentry{Illinois} 
            \addlegendentry{South Carolina} 
            \addlegendentry{Texas} 

        \end{axis}
    \end{tikzpicture}
\end{center}
    \caption{Percentage of hourly energy requirement of EVs. }\label{fig:drivingprofiles}
\end{figure}
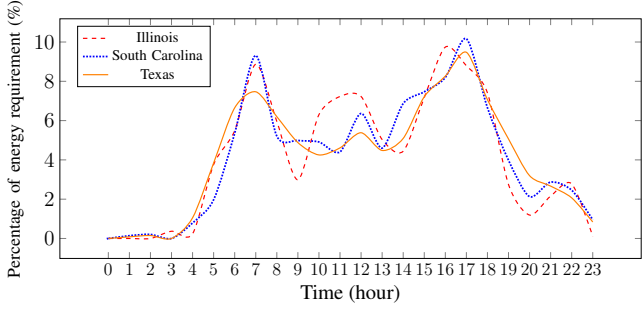

Since the data is only available for one typical day, we use the same EV driving profiles for both summer and winter, ignoring possible seasonal variations in EV energy consumption. Fig. \ref{fig:drivingprofiles} provides an overview of these driving profiles. We will explain how to integrate these driving profiles into our model in detail in 
Section~\ref{Experimental Setup}.
\section{Computational Experiments} \label{sec:MOPFCompExperiments} 
This section presents the results of our computational experiments, which have been conducted to investigate the effectiveness of our approach and assess the aggregate impact of EV charging load on optimal network operations under different scenarios. We first discuss the experimental setup to explain how we integrate input datasets in Section \ref{sec:MOPFInput} into the MOPF model in Section \ref{sec:MINLP}. Then, we graphically illustrate the results of the computational experiments.

\subsection{Experimental Setup} \label{Experimental Setup}

This section presents how we process input datasets so as to be compatible with the multi-period formulation. For the EV-aware MOPF model, we consider a 24-hour period, in 1-hour slots, from 00:00 to 00:00 the following day.

The TAMU instances include network parameters for only a single period. Most of the parameters remain constant throughout the day, but some parameters vary on an hourly basis, such as power load and electricity price. To extend to the multi-period version, we multiply the real and reactive power load of the grid with the corresponding normalized values as follows:
$$
 p_{it}^d = p_i^d \times  \overline{u}_t, \  q_{it}^d = q_i^d \times  \overline{u}_t \quad i \in \mathcal{B}, t \in \mathcal{T},
 $$
$$
  p_{it}^d = p_i^d \times  \overline{w}_t, \  q_{it}^d = q_i^d \times  \overline{w}_t \quad i \in \mathcal{B}, t \in \mathcal{T}.
$$
Here $\overline{u}_t$ and $\overline{w}_t$ are the normalized values for summer and winter seasons (see Fig. \ref{fig:Demand}).
Since hourly electricity price data is not publicly available, we assume that the price is constant over the 24-hour period 
\cite{eia.gov2020}. We adjust emission parameters $e_{t}$ to corresponding marginal emission factors in Section \ref{sec:marginalemission}.



We utilize EV driving profiles to adjust parameter settings related to EV charging.  We assume that EVs can only be charged at a bus with a nonzero real power load. If there is no real load at bus $i$, we set $c_{it}=\overline{a}_{it}=\overline{b}_{it}=\underline{s}_{it}=\overline{s}_{it}=0 $ for all $ t \in \mathcal{T}$. 

Instead of individual EVs, we consider EV groups, each containing EVs with the same driving profile. These EV groups are assigned to buses according to their weighted energy demand, $ \sum_{t \in \mathcal{T} } w_v \Omega_{vt}$, in such a way that each bus $i \in  \mathcal{D}$ has one EV group. Here, $w_v$ represents a weighting factor obtained from the NHTS dataset. Assuming that if the bus load is high, the EV demand around the bus is also high, we assign an EV group with a higher demand to a bus with a higher load.


Now, each bus with a real power load has one group of EV, and thus we can use bus indices $i$ instead of vehicle indices $v$ for brevity. We adjust parameter settings for EV groups instead of individual EVs. We first set energy requirements $c_{it}$ to $w_i \Omega_{it}$ for $i \in \mathcal{D}, t \in \mathcal{T}$. According to the energy requirements, we determine stock parameter $\underline{s}_{it}$ for each period. If an EV has demand during successive periods $(t,t+1,\dots,t+n)$, the amount of battery charge at the beginning of the period $t$ should be sufficient to satisfy the demand in these successive periods.



As previously mentioned in Section \ref{sec:Drivingprofiles}, we use the Nissan Leaf as a baseline vehicle. The Nissan Leaf has a usable battery capacity of 32kWh with a 6.6kW onboard charger \cite{nissanleaf}. We assume that charging rates are constant over a period. For each bus $i \in \mathcal{D}$, the maximum limit of charging and discharging powers, $\overline{a}_{it}$ and $\overline{b}_{it}$, is set to $w_i\times6.6$. We set the stock parameter $\overline{s}_{it}$ to usable battery capacity, that is $w_i\times32$, and the charging efficiency to 90\% \cite{richardson2011optimal}. 
We also set the initial battery SOC $I_i$ for EV at bus $i \in \mathcal{D}$ to zero.





To make the ratio of the electricity consumption and EV charging demand consistent with the actual data, we calculate a parameter called \texttt{weight} as the ratio between the total power demand of the grid and daily electricity consumption as follows:
\begin{equation}\label{weight}
\texttt{weight} = \frac{ \sum_{i \in \mathcal{B}} \sum_{t \in \mathcal{T}} p_{it}^d  }{\max(\sum_{t \in \mathcal{T}} u_t, \sum_{t \in \mathcal{T}} w_t)}.
\end{equation}
Here, $u_t$ and $w_t$ are the electricity consumption for summer and winter obtained by normalization. 
We multiply EV charging requirement $c_{it}$ and thereby all charging parameters $\overline{a}_{it}, \overline{b}_{it},\underline{s}_{it}, \overline{s}_{it},I_i$ with \texttt{weight}. Note that we use the same \texttt{weight} for different seasons. Table \ref{tab:demand} summaries the daily demand data at a macro level.

\begin{table}[ht]
\centering
\caption{Daily demands in kWh.}\label{tab:demand}
\begin{tabular}{|l|rrr|}
\hline
& \multicolumn{3}{|c|}{Region}  \\
&  IL &  SC &  TX  \\ \hline 
August demand ($\sum_{t \in \mathcal{T}} u_t$)
& 149690 & 82575 & 1256288  \\ 
December demand ($\sum_{t \in \mathcal{T}} w_t$)
& 129297 & 69901 & 923604  \\
Grid demand ($\sum_{i \in \mathcal{B}}  p_{i}^d$)
& 1476 & 7751 & 67109  \\
Total grid demand ($\sum_{i \in \mathcal{B}} \sum_{t \in \mathcal{T}} p_{it}^d $)
& 29276 & 151214 & 1313994  \\
EV demand ($\sum_{t \in \mathcal{T}} \sum_{v \in \mathcal{V}} \Omega_{vt}w_v$)
& 9632 & 3312 & 22646  \\  \hline 
\texttt{weight} & 0.19 & 1.83 & 1.05  \\
\hline 
\end{tabular}
\end{table}

We assume that an EV cannot be driven and connected to the grid in the same period. If the EV is in driving mode at time $t$, we set $a_{it}$ and $b_{it}$ to be zero. If EV has no energy requirement at time $t$, we assume that it is connected to the grid, and optimization decides exactly one mode among the following three based on the driving profile: charging, discharging, or rest. It can only be in one of these modes during time~$t$. In addition, if the operating cost of a generator is zero, then we set the lower bound on the power generation as $\underline{p}_i=0$.

The proposed MOPF model allows bidirectional flow, G2V, and V2G. In the case where only G2V is allowed, the maximum allowable discharging rate $\overline{b}_{it}$  is set to zero.

\subsection{Results}

In this section, we report the result of our computational experiments on three instances from PGLIB-OPF. A 64-bit desktop with Intel Xeon CPU E5-2630 v4 with a 2.20GHz processor and 32 GB RAM is used for all experiments. Our code is written in Python programming language using Spyder environment. The solvers GUROBI and IPOPT are used to solve the SOCP relaxation and the NLP models, respectively. 

We perform Algorithm \ref{Algorithm1} for each region, considering four different cases depending on the season of the year and direction of flow, namely G2V in summer, G2V-V2G in summer, G2V in winter, G2V-V2G in winter. We also create a benchmark strategy for comparison purposes, where we assume all EV drivers start charging their EVs at midnight. 

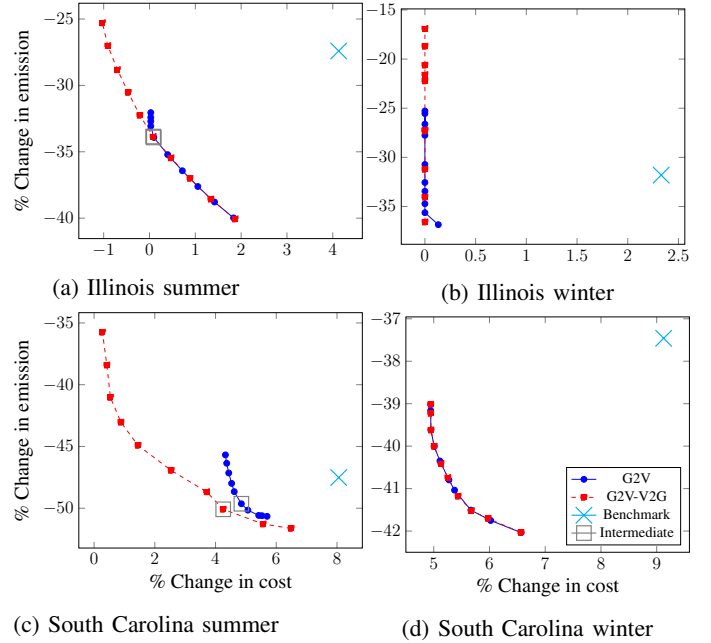
\begin{figure}[ht]
\begin{subfigure}{0.20\textwidth}
  \centering
  \begin{tikzpicture}[scale=0.55]
\begin{axis}[
    ylabel={\Large \% Change in emission},
    y tick label style={font=\large},
    x tick label style={font=\large},
    legend pos=north west,
    grid style=dashed,
]

\addplot[
    color=blue,
    mark=*,
    ]
    coordinates {
   (1.8348371423864926,-39.981623344996216)
(1.4188381381675004,-38.782709778344426)
(1.0548754724642702,-37.615584293176)
(0.7169321387720609,-36.42935284885862)
(0.395974781135622,-35.21397835405965)
(0.09478086123281036,-33.92201005171948)
(0.028248680924946735,-33.07275064266871)
(0.02824817313891446,-32.7066520319374)
(0.028261765157982954,-32.41544526125858)
(0.028249540230728502,-32.0500596362918)
    };
    
    \addplot[
    color=gray,
    mark=square,
    mark size=5pt,
    thick
    ]
    coordinates {
(0.07900713678668128,-33.85162483119248)

    };
    \addplot[
    color=gray,
    mark=square,
    mark size=5pt,
    thick
    ]
    coordinates {
(0.09478086123281036,-33.92201005171948)

    };
\addplot[
    color=red,
    mark=square*,
    style=dashed,
    ]
    coordinates {
   (1.8681992835998176,-40.06231314205295)
(1.3383867866801265,-38.56403202942513)
(0.8840807649302028,-37.00482842857282)
(0.4682620317608318,-35.45593037540387)
(0.07900713678668128,-33.85162483119248)
(-0.21147121182275966,-32.253085123467486)
(-0.4632644693176521,-30.532928613034557)
(-0.6998018237607389,-28.81901854152922)
(-0.9103219888416689,-27.00516424408736)
(-1.0269947768701775,-25.288108091614568)
    };

 \addplot[
    color=cyan,
    mark=x,
    mark size=8pt,
    thick
    ]
    coordinates {
(4.129300209307533,-27.407726035179838)

    };\end{axis}

\end{tikzpicture}
\caption{Illinois summer }\label{fig:IL-SCPareto_a}
\end{subfigure}
\begin{subfigure}{.35\textwidth}
  \centering
 
  \begin{tikzpicture}[scale=0.55, transform shape]
\begin{axis}[
    y tick label style={font=\large},
    x tick label style={font=\large},
    legend pos=north west,
    ymajorgrids=false,
    grid style=dashed,
]

\addplot[
    color=blue,
    mark=*,
    ]
    coordinates {
    (0.1320554253347133,-36.8500615942225)
(1.8587332913382023e-09,-35.62934831136246)
(2.2836791213817774e-09,-34.71974016366886)
(2.2836791213817774e-09,-33.45007747746776)
(1.862172112829729e-09,-32.55716577227321)
(2.2871380529405057e-09,-30.720306832167502)
(2.2871380529405057e-09,-27.757311972407855)
(2.2871380529405057e-09,-25.280661476072535)
(2.2871380529405057e-09,-26.623824167766628)
(2.2871380529405057e-09,-25.547606766310153)
    };
    
\addplot[
    color=red,
    mark=square*,
    style=dashed,
    ]
    coordinates {
   (0.0012316014995455622,-36.57995907184723)
(2.2802201898230487e-09,-34.05105515392991)
(2.2836791213817774e-09,-31.222539737727377)
(2.2836791213817774e-09,-27.320870112609505)
(2.2836791213817774e-09,-27.241063132856876)
(2.2836791213817774e-09,-22.21970135033301)
(2.2836791213817774e-09,-21.57869965557575)
(2.2871380529405057e-09,-20.613532458754616)
(2.2871380529405057e-09,-18.698642060603696)
(2.2871380529405057e-09,-16.890454794944205)
    };
 \addplot[
    color=cyan,
    mark=x,
    mark size=8pt,
    thick
    ]
    coordinates {
(2.3298033100700026,-31.80585396077706)

    };
    

\end{axis}
\end{tikzpicture} 
  \caption{ Illinois winter }\label{fig:IL-SCPareto_b}
\end{subfigure}

\begin{subfigure}{0.20\textwidth}
  \centering
 
  \begin{tikzpicture}[scale=0.55]
\begin{axis}[
    xlabel={ \Large \% Change in cost},
    ylabel={ \Large \% Change in emission},
    y tick label style={font=\large},
    x tick label style={font=\large},
    legend pos=north west,
    grid style=dashed,
]

\addplot[
    color=blue,
    mark=*,
    ]
    coordinates {
    (5.696868373545408,-50.6558904148223)
(5.5297525623212005,-50.60573353505383)
(5.422577435420465,-50.58298080761577)
(5.065193351091259,-50.153982696254666)
(4.8508042606385,-49.63867422143416)
(4.6151643180063795,-48.65753141016942)
(4.52978993650081,-47.989752186507324)
(4.434409670683222,-47.145376678772)
(4.367578058961336,-46.36026604602246)
(4.325095087556578,-45.68282976009954)
    };

\addplot[
    color=red,
    mark=square*,
    style=dashed,
    ]
    coordinates {
(6.489632143714304,-51.63595005195853)
(5.568373756533885,-51.278249081290696)
(4.252524219854367,-50.08761279411309)
(3.704201925005115,-48.66873758838299)
(2.539398015746234,-46.94053574490443)
(1.4480599383499893,-44.91493610774872)
(0.8825899678979935,-43.038323612086145)
(0.5397464537843446,-41.03156695155724)
(0.4197464537843446,-38.40718760495678)
(0.27478622680040715,-35.75704916681376)
    };    
     \addplot[
    color=cyan,
    mark=x,
    mark size=8pt,
    thick
    ]
    coordinates {
(8.057319849922711, -47.51787329684767)

    };
    
    \addplot[
    color=gray,
    mark=square,
    mark size=5pt,
    thick
    ]
    coordinates {
(4.8508042606385,-49.63867422143416)
};

    \addplot[
    color=gray,
    mark=square,
    mark size=5pt,
    thick
    ]
    coordinates {
(4.252524219854367,-50.08761279411309)
};


\end{axis}
\end{tikzpicture}

\caption{South Carolina summer } \label{fig:IL-SCPareto_c} 
\end{subfigure}
\begin{subfigure}{.35\textwidth}
  \centering
  \begin{tikzpicture}[scale=0.55, transform shape]
\begin{axis}[
    xlabel={\Large \% Change in cost},
     y tick label style={font=\large},
    x tick label style={font=\large},
    xmin=4.50,
    ymin=-42.5, 
    legend pos=south east,
    ymajorgrids=false,
    grid style=dashed,
]

\addplot[
    color=blue,
    mark=*,
    ]
    coordinates {
   (6.569393243582627,-42.03659878973923)
(6.021191395429295,-41.751904799546615)
(5.666617359724065,-41.507981998384714)
(5.37334286214752,-41.036731504571996)
(5.272095594458021,-40.79512686582607)
(5.110801759714496,-40.34960394739813)
(5.008867439716306,-40.00948680041029)
(4.949460146013213,-39.61743642458768)
(4.945158055583467,-39.162449610998365)
(4.944967475241416,-39.020245628821435)
    };
    
\addplot[
    color=red,
    mark=square*,
    style=dashed,
    ]
    coordinates {
      (6.563528512372708,-42.03375011783914)
(5.975413246261896,-41.704204620269735)
(5.675235120318561,-41.51918299050861)
(5.437634454874512,-41.18049574814257)
(5.253066050840572,-40.73880601009412)
(5.1297644942124885,-40.409236355893704)
(5.007978654884666,-40.00516337869629)
(4.949533151179254,-39.61916080637172)
(4.944995513339557,-39.228353380943595)
(4.944978997519207,-39.01514006502542)
    };
    
    \addplot[
    color=cyan,
    mark=x,
    mark size=8pt,
    thick
    ]
    coordinates {
(9.12475057822829, -37.45978890193726)

    };
     \addplot[
    color=gray,
    mark=square,
    mark size=5pt,
    thick
    ]
    coordinates {
(0,-50)

    };
       
    \addlegendentry{G2V} 
    \addlegendentry{G2V-V2G}
    \addlegendentry{Benchmark}
    \addlegendentry{Intermediate}
                
\end{axis}
\end{tikzpicture} 
 \caption{ South Carolina winter}\label{fig:IL-SCPareto_d}
\end{subfigure}

\caption{Pareto frontier of \% changes in costs and emissions.}\label{fig:IL-SCPareto}
\end{figure}
\begin{figure*}[ht]
\allowdisplaybreaks
\begin{subfigure}{0.33\textwidth}

\begin{tikzpicture}[scale=0.45, transform shape]
        \begin{axis}[
            xlabel = \Large Time (hour),
            ylabel = \Large Energy (kWh) (normalized) ,
            title style={font=\large},
            xtick = {0,1,2,3,4,5,6,7,8,9,10,11,12,13,14,15,16,17,18,19,20,21,22,23},
            x tick label style={font=\large},
            y tick label style={font=\Large},
            ymin = 0,
            ymax = 1.1,
            x = 0.42cm,
            legend pos =south east,
            xmajorgrids=false,
            grid style=dashed]
            
        \addplot[smooth, thick, color = blue, mark = *,style=dashed]
            plot coordinates {
(0,0.769417373611305)
(1,0.7182193288786001)
(2,0.6823067537509248)
(3,0.6567355251706724)
(4,0.6373911448165036)
(5,0.6341376203815321)
(6,0.6401673894500886)
(7,0.6658353421553957)
(8,0.6799355997839838)
(9,0.7200498874909985)
(10,0.7802331819754933)
(11,0.8397632018713678)
(12,0.8983259875561163)
(13,0.9451649856458716)
(14,0.9776457077224721)
(15,1.0008155257299371)
(16,1.014762860665576)
(17,1.0206078975911492)
(18,1.020974184099852)
(19,0.9942329656398512)
(20,0.9600391705453402)
(21,0.9359179484208648)
(22,0.8961759062230155)
(23,0.834422388653523)
            };
        
        \addplot[smooth, thick, color = purple, mark = *]
            plot coordinates {
                (0,0.007359949178878409)
(1,0.058560399842377264)
(2,0.09447822455328324)
(3,0.14663753654471254)
(4,0.13947676217210908)
(5,0.12926575754888323)
(6,0.1010357103551724)
(7,0.09114164798331546)
(8,0.04198441790463244)
(9,0.014522553683687845)
(10,0.0012448088288639228)
(11,0.006126996108116823)
(12,0.0042523242192927775)
(13,0.005558216344773747)
(14,0.0031391095272924234)
(15,3.03171324851656e-06)
(16,0.0026381178394796136)
(17,0.00045487553439062717)
(18,0.0014327600170663934)
(19,0.002062254345622519)
(20,0.0032494399150685507)
(21,6.048999972266457e-06)
(22,0.0016935444118607113)
(23,1.027064041016553e-06)
                
            };

            \addplot[smooth, thick, color = orange, mark = *,style=dotted]
            plot coordinates {
    (23,100)
            };
            
        \end{axis}
    \end{tikzpicture}
\caption{Cost minimization G2V}
  \label{fig:HourlyLoadCostMin}
\end{subfigure}
\begin{subfigure}{0.33\textwidth}
\begin{tikzpicture}[scale=0.45, transform shape]
        \begin{axis}[
            xlabel = \Large Time (hour),
            title style={font=\large},
            xtick = {0,1,2,3,4,5,6,7,8,9,10,11,12,13,14,15,16,17,18,19,20,21,22,23},
            x tick label style={font=\large},
            ytick = {0, 0.2, 0.4, 0.6, 0.8, 1,1.2},
            y tick label style={font=\Large},
            ymin = 0,
            x = 0.42cm,
            legend pos = south east,
            xmajorgrids=false,
            grid style=dashed]
        \addplot[smooth, thick, color = purple, mark = *]
            plot coordinates {
                (0,0.007433639222523393)
(1,0.05874892897097661)
(2,0.09467876423816221)
(3,0.20742075436547652)
(4,0.19073000756956385)
(5,0.08966928250055618)
(6,0.08224704160020158)
(7,0.05719436703398412)
(8,0.026173491872295077)
(9,0.010122704456278738)
(10,0.0011854833748108134)
(11,0.0034192997293434123)
(12,0.007069681445627979)
(13,0.005508131527286677)
(14,0.0031839036852220684)
(15,1.2835303244454717e-05)
(16,0.0026423528596021205)
(17,0.0004618719174190801)
(18,0.0014396559775360731)
(19,0.0020605897310680995)
(20,0.0032560519000698965)
(21,2.694677777368229e-06)
(22,0.0016957995556437185)
(23,3.648564020389501e-06)
                
            };
            \addplot[smooth, thick, color = blue, mark = *,style=dashed]
            plot coordinates {
                (0,0.7694031711226383)
(1,0.7181280526138084)
(2,0.6822247223320125)
(3,0.6583974098635963)
(4,0.6384161610742303)
(5,0.6327596084335357)
(6,0.6395586257739435)
(7,0.6646385779672146)
(8,0.6796157348691921)
(9,0.7199453338478544)
(10,0.7802316721604198)
(11,0.8396620435105686)
(12,0.8983142873340269)
(13,0.9451645094420589)
(14,0.9776462314684503)
(15,1.0008156333152114)
(16,1.0147628921083245)
(17,1.0206079156055499)
(18,1.0209737280235696)
(19,0.9942332647761473)
(20,0.9600392050385043)
(21,0.9359179820402672)
(22,0.8961758956899556)
(23,0.8344224583378884)
            };
        \end{axis}
    \end{tikzpicture}
    \caption{Intermediate case G2V}
  \label{fig:HourlyLoadInt}
\end{subfigure}
\begin{subfigure}{0.33\textwidth}
  
    \begin{tikzpicture}[scale=0.45, transform shape]
        \begin{axis}[
            xlabel = \Large Time (hour),
            title style={font=\large},
            xtick = {0,1,2,3,4,5,6,7,8,9,10,11,12,13,14,15,16,17,18,19,20,21,22,23},
            x tick label style={font=\large},
            ytick = {0, 0.2, 0.4, 0.6, 0.8, 1,1.2},
            y tick label style={font=\Large},
            ymin = 0,
            x = 0.42cm,
            legend pos = south east,
            xmajorgrids=false,
            grid style=dashed]
        \addplot[smooth, thick, color = purple, mark = *]
            plot coordinates {
                (0,0.007717371778680091)
(1,0.1112030252013778)
(2,0.1788946755212732)
(3,0.24083742681329012)
(4,0.20351131529818994)
(5,0.04918595300321062)
(6,0.000628389679106025)
(7,0.01544881272720977)
(8,5.342141100379281e-06)
(9,0.00034576045446333575)
(10,3.862565177523163e-06)
(11,0.0010193024220560118)
(12,0.013371136404437515)
(13,0.012054839636279553)
(14,0.005369629010589278)
(15,0.000727745829433268)
(16,0.006649025015137076)
(17,0.000482232078570991)
(18,0.006198417463610842)
(19,0.0006074449591976392)
(20,0.002063031875773492)
(21,8.048783091170444e-07)
(22,1.7879053249800037e-06)
(23,7.955777996659286e-07)
                
            };
            \addplot[smooth, thick, color = blue, mark = *,style=dashed]
            plot coordinates {
                (0,0.7694938200185634)
(1,0.719459201700468)
(2,0.6843672628806192)
(3,0.6592090371461561)
(4,0.6388798271563234)
(5,0.6320300562403721)
(6,0.6388465667068839)
(7,0.6644838218961043)
(8,0.6798118734793065)
(9,0.7197475639158286)
(10,0.7802063578610171)
(11,0.8395833462640032)
(12,0.8982681360278206)
(13,0.9451914690697745)
(14,0.9776685788803874)
(15,1.0008219959709768)
(16,1.0147473529378466)
(17,1.0206078288472897)
(18,1.0209158284447057)
(19,0.9943119003639003)
(20,0.9600280136883568)
(21,0.9359179860569936)
(22,0.8961423848673339)
(23,0.8344223825337824)
                
            };
            
        \end{axis}
    \end{tikzpicture}
    \caption{Emission minimization G2V}
  \label{fig:HourlyLoadEmisMin}
\end{subfigure}

\begin{subfigure}{0.33\textwidth}
 
\begin{tikzpicture}[scale=0.45, transform shape]
        \begin{axis}[
            xlabel = \Large Time (hour),
            ylabel = \Large Energy (kWh) (normalized) ,
            title style={font=\large},
            xtick = {0,1,2,3,4,5,6,7,8,9,10,11,12,13,14,15,16,17,18,19,20,21,22,23},
            x tick label style={font=\scriptsize},
            ymin = 0,
            ymax = 1.1,
            x = 0.42cm,
            legend pos = south east,
            xmajorgrids=false,
            grid style=dashed]
        \addplot[smooth, thick, color = purple, mark = *]
            plot coordinates {
(0,0.056029873247570694)
(1,0.12691289577882622)
(2,0.1631391236623838)
(3,0.200109276412299)
(4,0.20609828957818782)
(5,0.17785011013888732)
(6,0.1619854707150826)
(7,0.14428955974193647)
(8,0.11183779184788815)
(9,0.08536569385231133)
(10,0.05482563658335404)
(11,0.03385213508859682)
(12,0.020878818671510545)
(13,0.0052079614662105846)
(14,0.00210126448485158)
(15,4.500231192186232e-05)
(16,0.0016081472631366224)
(17,0.0008664656713155605)
(18,0.001293390367510455)
(19,0.002323422062976807)
(20,0.0036159321052318285)
(21,5.056249788367087e-05)
(22,0.004676035158959557)
(23,3.944159344087082e-07)
                
            };
            \addplot[smooth, thick, color = orange, mark = *,style=dotted]
            plot coordinates {
(0,0.0)
(1,0.0)
(2,0.0)
(3,0.0)
(4,0.0)
(5,0.0)
(6,0.0)
(7,0.0)
(8,0.0)
(9,0.0)
(10,0.0)
(11,0.0)
(12,0.00047240657168890626)
(13,0.0306140548831678)
(14,0.06219769331912995)
(15,0.08513620115811701)
(16,0.10093707111645385)
(17,0.10675486939250971)
(18,0.10735771904227195)
(19,0.0810481722339824)
(20,0.045973324737893353)
(21,0.017228046919143872)
(22,4.6623942336063947e-07)
(23,2.716784232522409e-05)
                
            };
            \addplot[smooth, thick, color = blue, mark = *,style=dashed]
            plot coordinates {
(0,0.7717682318838535)
(1,0.7203643741575818)
(2,0.6843197490031749)
(3,0.6583267426778034)
(4,0.6397520663328335)
(5,0.6366748249102886)
(6,0.6425218857998358)
(7,0.668752766032579)
(8,0.6823459599380362)
(9,0.7229695192696267)
(10,0.7825404968649116)
(11,0.8404967023977048)
(12,0.8979020105132896)
(13,0.9175455974273021)
(14,0.9214562009973699)
(15,0.9236418490064269)
(16,0.9233720779169744)
(17,0.9240104905926226)
(18,0.9237417510391052)
(19,0.9208145457557104)
(20,0.9185730886934027)
(21,0.9203695951482949)
(22,0.8962101709971161)
(23,0.8344004424755622)
            };
            
        \end{axis}
    \end{tikzpicture}
    \caption{Cost minimization G2V and V2G}
  \label{fig:HourlyLoadCostMin_v2g}
\end{subfigure}
\begin{subfigure}{0.33\textwidth}
 
\begin{tikzpicture}[scale=0.45, transform shape]
        \begin{axis}[
            xlabel = \Large Time (hour),
            title style={font=\large},
            xtick = {0,1,2,3,4,5,6,7,8,9,10,11,12,13,14,15,16,17,18,19,20,21,22,23},
            x tick label style={font=\large},
            ytick = {0, 0.2, 0.4, 0.6, 0.8, 1,1.2},
            y tick label style={font=\Large},
            ymin = 0,
            x = 0.42cm,
            legend pos = south east,
            xmajorgrids=false,
            grid style=dashed]
        \addplot[smooth, thick, color = blue, mark = *, style=dashed]
            plot coordinates {
                (0,0.7710366571606718)
(1,0.722000153808532)
(2,0.6850519027839357)
(3,0.6591639908465855)
(4,0.6400466599859925)
(5,0.6428162989048755)
(6,0.6345011823734849)
(7,0.6612994309656067)
(8,0.6752672324971437)
(9,0.693938926928186)
(10,0.7753719876805067)
(11,0.8010273679012138)
(12,0.892284169344823)
(13,0.9322283492022511)
(14,0.9351053385524196)
(15,0.9226433111824358)
(16,0.9590973556583405)
(17,0.9225328393883295)
(18,0.9651003673130751)
(19,0.9198120296392903)
(20,0.9177880697453654)
(21,0.9194884514691485)
(22,0.8932728969426906)
(23,0.8344019048222517)
                
            };
            \addplot[smooth, thick, color = purple, mark = *,style=solid]
            plot coordinates {
               (0,0.08705607307965067)
(1,0.19253476005467313)
(2,0.22814148510350296)
(3,0.24201090876742584)
(4,0.2342851752121926)
(5,0.14184765880909486)
(6,0.10416645124899397)
(7,0.10789945360982298)
(8,0.0423806120463901)
(9,0.05097862349495039)
(10,0.006062857007852545)
(11,0.04220917150399759)
(12,0.02274447127576313)
(13,0.005368472642250401)
(14,0.0027501161920532456)
(15,0.0007150864014623169)
(16,0.0024556382591582595)
(17,0.0008101680652820733)
(18,0.0018497377242184697)
(19,0.002118482870165034)
(20,0.0031432795631278557)
(21,4.81045628438564e-05)
(22,0.004618529540003422)
(23,4.981566065230309e-09)
            };
            
\addplot[smooth, thick, color = orange, mark = *,style=dotted]
            plot coordinates {
               (0,0.0)
(1,0.0)
(2,0.0)
(3,0.0)
(4,0.0)
(5,0.0)
(6,0.0)
(7,0.0)
(8,1.6188409001285916e-07)
(9,7.156812143091053e-08)
(10,6.24244377711174e-07)
(11,2.9197510316387697e-07)
(12,4.383415046660865e-06)
(13,0.02779758725717735)
(14,0.060399780544300015)
(15,0.08661257233229705)
(16,0.08507747478618932)
(17,0.10872573879249782)
(18,0.08683745498771878)
(19,0.0822152522255539)
(20,0.0470063972221821)
(21,0.01814778902288754)
(22,7.291568084250182e-07)
(23,2.3903911409504448e-05)
            };
        \end{axis}
    \end{tikzpicture}
    \caption{Intermediate case G2V and V2G}
  \label{fig:HourlyLoadCostMinInt_v2g}
\end{subfigure}
\begin{subfigure}{0.33\textwidth}
 
    \begin{tikzpicture}[scale=0.45, transform shape]
        \begin{axis}[
            xlabel = \Large Time (hour),
            title style={font=\large},
            xtick = {0,1,2,3,4,5,6,7,8,9,10,11,12,13,14,15,16,17,18,19,20,21,22,23},
            x tick label style={font=\large},
            ytick = {0, 0.2, 0.4, 0.6, 0.8, 1,1.2},
            y tick label style={font=\Large},
            ymin = 0,
            x = 0.42cm,
            legend pos = south east,
            legend style={at={(0.98,0.3)}},
            xmajorgrids=false,
            grid style=dashed]
        \addplot[smooth, thick, color = blue, mark = *,style=dashed]
            plot coordinates {
                (0,0.7695038980320936)
(1,0.7209701853179434)
(2,0.6851351237346829)
(3,0.659215184385958)
(4,0.6400285280011223)
(5,0.6324206373561863)
(6,0.6388393850436804)
(7,0.6644822977049682)
(8,0.6793073307216458)
(9,0.7197311108303235)
(10,0.7074083732367854)
(11,0.8395468715665796)
(12,0.8981786555976102)
(13,0.945188498577395)
(14,0.977667310676502)
(15,1.0008131226396801)
(16,1.0147534500413262)
(17,1.0205985858394817)
(18,1.0209130611569537)
(19,0.9938466324953912)
(20,0.9600126908757952)
(21,0.8595713927461252)
(22,0.8961121909072529)
(23,0.8343891889177883)
                
            };
            \addplot[smooth, thick, color = purple, mark = *,style=solid]
            plot coordinates {
                (0,0.007171802642824964)
(1,0.20285999154177844)
(2,0.2366447448989197)
(3,0.24201910187162376)
(4,0.23430260169062034)
(5,0.048535305552679466)
(6,7.52772102258941e-05)
(7,0.012825335554298987)
(8,4.56408419211983e-06)
(9,0.0003297576364144277)
(10,1.0802061586796723e-06)
(11,0.0010200425402901389)
(12,0.017462584401553873)
(13,0.01199437790361371)
(14,0.005331974820328002)
(15,0.0007276720562983801)
(16,0.007189416826117191)
(17,0.00048488188620378787)
(18,0.006201574396366165)
(19,0.0006065401902036518)
(20,0.002066025793782603)
(21,0.0)
(22,1.663579668241631e-05)
(23,0.0)
                
            };
            
\addplot[smooth, thick, color = orange, mark = *,style=dotted]
            plot coordinates {
                (0,0.0)
(1,0.0)
(2,0.0)
(3,0.0)
(4,0.0)
(5,4.067938624411376e-08)
(6,1.2257540188840405e-06)
(7,1.0509451838090455e-07)
(8,0.0005669697717035076)
(9,1.908108957343184e-05)
(10,0.07837999470994837)
(11,4.101986832020259e-05)
(12,0.0)
(13,1.029421845633577e-06)
(14,2.0241415892095003e-06)
(15,1.0888414024094766e-05)
(16,1.3134377984617485e-06)
(17,1.1370837830429164e-05)
(18,4.247267138161795e-06)
(19,0.0005232867394915296)
(20,1.8889872037706404e-05)
(21,0.08371593039101242)
(22,3.574492335824073e-05)
(23,3.7852618641211795e-05)
                
            };
            
        \addlegendentry{Electricity generation excluding EVs}
          \addlegendentry{Electricity generated for charging EVs}
        \addlegendentry{V2G power profile} 
        
        \end{axis}
    \end{tikzpicture}
    \caption{Emission minimization G2V and V2G}
  \label{fig:HourlyLoadEmisMin_v2g}
\end{subfigure}
    \caption{Hourly load variations for South Carolina Summer (Intermediate cases are indicated in Fig. \ref{fig:IL-SCPareto_c}). }
    \label{fig:HourlyLoadSC}
\end{figure*}
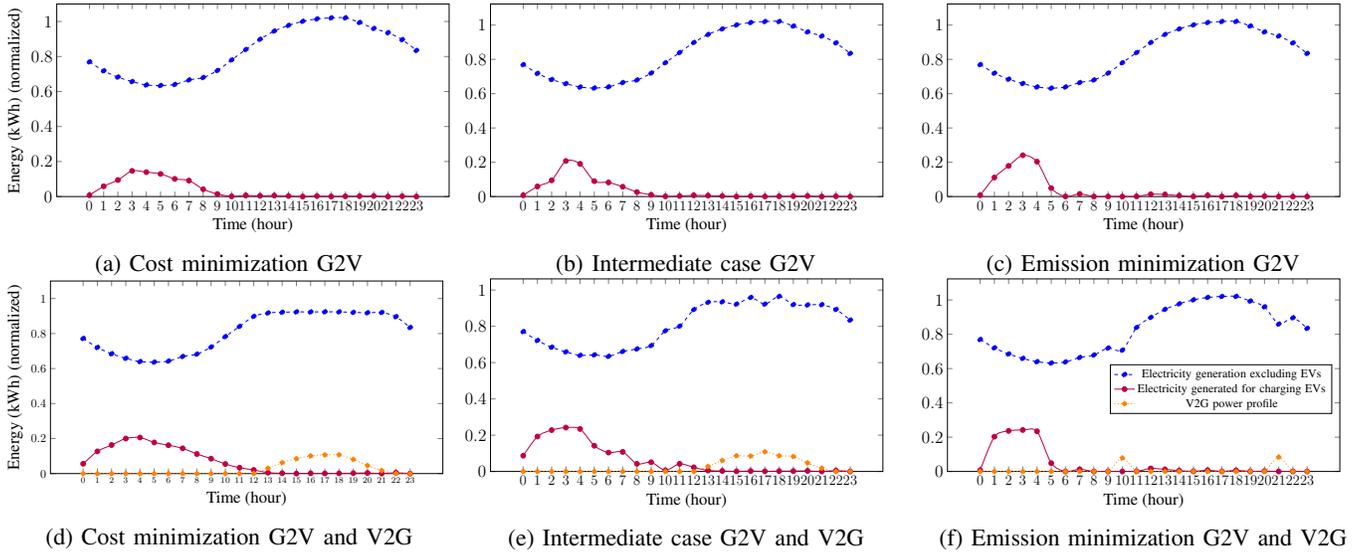

We illustrate how $\%$ changes in total generation cost and marginal emission move in relation to each other in Fig. \ref{fig:IL-SCPareto}. This figure shows the G2V case (in solid blue line) and G2V-V2G case (in red dashed line) computed via NLP models. The cyan cross represents the corresponding result for the benchmark case.
We calculate the change in cost as  $\%\text{Change}=100\times(\text{generation cost with EVs}/\text{generation cost without EVs}-1)$. For $\%$ change in emission, we first assume all vehicles in our model belong to the average gasoline car class \cite{afdcfuel} and calculate their resulting emission. We then provide the changes in emission as 
$\%\text{Change}=100\times(\text{marginal emission by EVs}/\text{emission by gasoline cars}-1)$. In Figs. \ref{fig:HourlyLoadIL} and \ref{fig:HourlyLoadSC}, we plot the hourly electricity generation excluding EVs (in dashed blue line), hourly electricity generated for charging EVs (in solid red line), and hourly V2G power profile (in orange dotted line) for cost minimization, intermediate, and emission minimization cases.

\subsubsection{Cost-Emission Trade-off}
From Fig. \ref{fig:IL-SCPareto}, one can observe that changes in cost and emission are inversely correlated. When we set tight restrictions on emission, the cost tends to rise; while we aim to minimize cost, emission increases.  For example, as depicted in Fig. \ref{fig:IL-SCPareto_a}, emission decreases by approximately 3\% in Illinois during summer for the  G2V case with no change in the generation cost. Similar interpretations can also be made for the results observed in Figs. \ref{fig:IL-SCPareto_c} and \ref{fig:IL-SCPareto_d}. When we integrate the V2G  concept, we can decrease the generation cost much more, albeit with an increase in emission (see Figs. \ref{fig:IL-SCPareto_a} and \ref{fig:IL-SCPareto_c}).

In some cases, despite the tighter restrictions on emission, the cost remains the same up to a point. At this point, the reduction of emission is \textit{free}. This can be clearly observed in Fig. \ref{fig:IL-SCPareto_b}. 
In this case, we do not observe significant changes in cost since network demand is lower in winter than summer (see  Table \ref{tab:demand}), and generator lower bounds almost match the demand in winter. However, by coordinating the charging time of EVs, we can achieve emission savings.

For the 2000-bus Texas case, since the added EV demand compromise a smaller part (1.8\%) of the total network demand, different charging patterns do not significantly affect the solution in terms of cost and emission. In particular, the solution that provides the least cost is almost identical to the solution that yields the least emission. Hence, we only provide optimal and benchmark values for the Texas case in Table \ref{tab:resultstexas} rather than the Pareto frontier. We note that when we artificially increase EV penetration level in the Texas case, we obtain Pareto Frontiers similar to Illinois and South Carolina test cases.
 

\begin{table}[ht]
\begin{center}
\caption{Results for Texas.}
    \label{tab:resultstexas}
\begin{tabular}{ |p{0.14\linewidth}|p{0.14\linewidth}|p{0.14\linewidth}|p{0.14\linewidth}|p{0.14\linewidth}| } 
 \hline
 & \multicolumn{2}{|c|}{Summer} & \multicolumn{2}{|c|}{Winter}   \\
 \hline
 & G2V-V2G & Benchmark & G2V-V2G & Benchmark \\ 
 Cost & \ \ 2.49  & \ \ 2.57 &  \ \ 3.65 & \ \ 3.68  \\ 
 Emission &-39.88  & -37.13 & -42.15 & -41.86 \\ 
 \hline
\end{tabular}
\end{center}
\end{table}


Finally, we would like to also note that, in each case, there exist multiple  points on the Pareto frontier that are strictly better than the benchmark strategy. This shows the promise of coordinated charging  in terms of both cost and emission.

\begin{figure*}[ht]
\allowdisplaybreaks

\begin{subfigure}{0.33\textwidth}

\begin{tikzpicture}[scale=0.45, transform shape]
        \begin{axis}[
            xlabel = \Large Time (hour),
            ylabel = \Large Energy (kWh) (normalized) ,
            title style={font=\large},
            xtick = {0,1,2,3,4,5,6,7,8,9,10,11,12,13,14,15,16,17,18,19,20,21,22,23},
            x tick label style={font=\large},
            y tick label style={font=\Large},
            ymin = 0,
            ymax = 1.1,
            x = 0.42cm,
            legend pos =south east,
            xmajorgrids=false,
            grid style=dashed]
            
        \addplot[smooth, thick, color = blue, mark = *,style=dashed]
            plot coordinates {
(0,0.8182734028455284)
(1,0.7597913868563685)
(2,0.7135089881436315)
(3,0.6826774078590786)
(4,0.6604339986449864)
(5,0.6475549749322492)
(6,0.6490293333333333)
(7,0.6677280643631436)
(8,0.6928099844173442)
(9,0.7289503143631436)
(10,0.7726676382113821)
(11,0.8183091214769648)
(12,0.8628393009349593)
(13,0.9025308420430216)
(14,0.9371380642750676)
(15,0.9700150819647696)
(16,0.990968166800813)
(17,1.0028769839096674)
(18,1.0076583607977372)
(19,0.9978629245564838)
(20,0.9743343949728998)
(21,0.9434283400948019)
(22,0.9208189138206432)
(23,0.8784292597144552)
            };
        
        \addplot[smooth, thick, color = purple, mark = *]
            plot coordinates {
                (0,0.0319101066395664)
(1,0.07989706368563687)
(2,0.02914789193766938)
(3,0.10058322425474255)
(4,0.1731833753387534)
(5,0.1914343543360434)
(6,0.20004346070460702)
(7,0.17340539837398372)
(8,0.1539206714092141)
(9,0.12468371815718156)
(10,0.0818195054200542)
(11,0.03658767662601626)
(12,0.0010496272493224932)
(13,6.7941358401084e-06)
(14,0.0016100373509485098)
(15,0.0028558807723577232)
(16,0.000631374527100271)
(17,5.526756985094851e-07)
(18,2.6568174390243903e-06)
(19,1.2488852506775068e-06)
(20,0.002859256111111111)
(21,4.337403279132792e-08)
(22,6.517664688346884e-08)
(23,2.8903431029810295e-08)
                
            };

            \addplot[smooth, thick, color = orange, mark = *,style=dotted]
            plot coordinates {
    (23,100)
            };
            
        \end{axis}
    \end{tikzpicture}
\caption{Cost minimization G2V}
  \label{fig:HourlyLoadCostMinIL}
\end{subfigure}
\begin{subfigure}{0.33\textwidth}
\begin{tikzpicture}[scale=0.45, transform shape]
        \begin{axis}[
            xlabel = \Large Time (hour),
            title style={font=\large},
            xtick = {0,1,2,3,4,5,6,7,8,9,10,11,12,13,14,15,16,17,18,19,20,21,22,23},
            x tick label style={font=\large},
            ytick = {0, 0.2, 0.4, 0.6, 0.8, 1,1.2},
            y tick label style={font=\Large},
            ymin = 0,
            x = 0.42cm,
            legend pos = south east,
            xmajorgrids=false,
            grid style=dashed]
        \addplot[smooth, thick, color = purple, mark = *]
            plot coordinates {
(0,0.04601556395663957)
(1,0.017622095528455284)
(2,1.5283684891598915e-05)
(3,0.004337007913279132)
(4,0.11894503929539296)
(5,0.21643686178861787)
(6,0.21520473035230353)
(7,0.19647655081300813)
(8,0.17140647086720867)
(9,0.13553051761517615)
(10,0.09170496815718156)
(11,0.05063565853658537)
(12,0.10521865989159891)
(13,4.935625020325204e-05)
(14,0.0016645562262872628)
(15,1.293203211382114e-05)
(16,7.14425216802168e-06)
(17,6.369740867208672e-06)
(18,0.0006429653902439025)
(19,0.002870502554200542)
(20,0.0028812378658536583)
(21,4.579744234417344e-06)
(22,2.72377466802168e-06)
(23,1.1451604878048778e-06)

            };
            \addplot[smooth, thick, color = blue, mark = *,style=dashed]
            plot coordinates {
              (0,0.8175599109756099)
(1,0.7590907160569106)
(2,0.7131221095642005)
(3,0.6814222330081301)
(4,0.6610365623306234)
(5,0.6471582730352304)
(6,0.6484185162601627)
(7,0.6671214275067752)
(8,0.6921921592140922)
(9,0.7280806463414634)
(10,0.7719305819783198)
(11,0.8175018376016261)
(12,0.8629591135501354)
(13,0.9025311919214092)
(14,0.937138627615176)
(15,0.9699062036762194)
(16,0.9909432290815311)
(17,1.002877063309844)
(18,1.0076845450284553)
(19,0.9979801190243902)
(20,0.9743345752913279)
(21,0.9434283676280624)
(22,0.9208189243693969)
(23,0.8784292637056645)
            };
        \end{axis}
    \end{tikzpicture}
    \caption{Intermediate case G2V}
  \label{fig:HourlyLoadIntIL}
\end{subfigure}
\begin{subfigure}{0.33\textwidth}
  
    \begin{tikzpicture}[scale=0.45, transform shape]
        \begin{axis}[
            xlabel = \Large Time (hour),
            title style={font=\large},
            xtick = {0,1,2,3,4,5,6,7,8,9,10,11,12,13,14,15,16,17,18,19,20,21,22,23},
            x tick label style={font=\large},
            ytick = {0, 0.2, 0.4, 0.6, 0.8, 1,1.2},
            y tick label style={font=\Large},
            ymin = 0,
            x = 0.42cm,
            legend pos = south east,
            xmajorgrids=false,
            grid style=dashed]
        \addplot[smooth, thick, color = purple, mark = *]
            plot coordinates {
(0,0.046071583468834686)
(1,0.010020663414634147)
(2,1.1569389905149051e-05)
(3,0.002848545711382114)
(4,0.09723407317073171)
(5,0.09883568699186991)
(6,0.21514588075880758)
(7,0.040296536111111114)
(8,0.05040834207317073)
(9,0.135546668699187)
(10,0.09427796815718158)
(11,0.20311866327913278)
(12,0.2464222723577236)
(13,0.011862844850948509)
(14,0.07773034891598916)
(15,0.010037764295392954)
(16,1.9551988482384825e-05)
(17,1.585439071815718e-05)
(18,0.006196688773712738)
(19,0.008537302574525745)
(20,0.0230452972899729)
(21,1.0528747086720868e-05)
(22,3.1407717005420057e-06)
(23,1.5369980352303523e-06)
                
            };
            \addplot[smooth, thick, color = blue, mark = *,style=dashed]
            plot coordinates {
                (0,0.8175842159891599)
(1,0.7586068555555555)
(2,0.7131222456507452)
(3,0.681443438705962)
(4,0.6591399424119242)
(5,0.6459998733062331)
(6,0.6484939024390245)
(7,0.6655795546747967)
(8,0.6911144404471544)
(9,0.7281291802168022)
(10,0.7719888875338753)
(11,0.8204876429539295)
(12,0.8679046714092141)
(13,0.9027525169376693)
(14,0.9386583062330623)
(15,0.9701584159214093)
(16,0.9909435225372629)
(17,1.002877410514431)
(18,1.0078471662398374)
(19,0.9981315172086721)
(20,0.9747920231707317)
(21,0.9434284868355691)
(22,0.9208189547567548)
(23,0.8784292766876015)
                
            };
            
        \end{axis}
    \end{tikzpicture}
    \caption{Emission minimization G2V}
  \label{fig:HourlyLoadEmisMinIL}
\end{subfigure}

\begin{subfigure}{0.33\textwidth}
 
\begin{tikzpicture}[scale=0.45, transform shape]
        \begin{axis}[
            xlabel = \Large Time (hour),
            ylabel = \Large Energy (kWh) (normalized) ,
            title style={font=\large},
            xtick = {0,1,2,3,4,5,6,7,8,9,10,11,12,13,14,15,16,17,18,19,20,21,22,23},
            x tick label style={font=\scriptsize},
            ymin = 0,
            ymax = 1.1,
            x = 0.42cm,
            legend pos = south east,
            xmajorgrids=false,
            grid style=dashed]
        \addplot[smooth, thick, color = purple, mark = *]
            plot coordinates {
(0,0.046024418699187)
(1,0.10419977981029811)
(2,0.14890964566395665)
(3,0.18030593360433608)
(4,0.2034165968834688)
(5,0.21644123780487803)
(6,0.21516453658536586)
(7,0.19648067547425474)
(8,0.17140971951219514)
(9,0.13550840040650408)
(10,0.0916710948509485)
(11,0.04610076788617886)
(12,0.0010524525474254743)
(13,3.7991548035230353e-06)
(14,0.0010228528658536586)
(15,0.001906700792682927)
(16,0.0008026860636856369)
(17,0.0006075216314363143)
(18,0.001239751951219512)
(19,0.0013396265040650406)
(20,0.0010422702235772356)
(21,3.638278319783198e-09)
(22,0.0)
(23,0.0)
                
            };
            \addplot[smooth, thick, color = orange, mark = *,style=dotted]
            plot coordinates {
(0,0.0)
(1,0.0)
(2,0.0)
(3,0.0)
(4,0.0)
(5,0.0)
(6,0.0)
(7,0.0)
(8,0.0)
(9,0.0)
(10,0.0)
(11,0.0)
(12,0.00012678387533875338)
(13,0.0001024159891598916)
(14,0.010103639566395664)
(15,0.03578521002710027)
(16,0.053276981707317074)
(17,0.06436360840108402)
(18,0.06985796747967479)
(19,0.06077636111111111)
(20,0.03894406978319783)
(21,0.013609224254742548)
(22,0.0010000731707317074)
(23,0.0003357960704607046)
                
            };
            \addplot[smooth, thick, color = blue, mark = *,style=dashed]
            plot coordinates {
            (0,0.8176333983739837)
(1,0.7594579464769649)
(2,0.7147492445799458)
(3,0.6833536368563685)
(4,0.6602433983739837)
(5,0.6472179471544715)
(6,0.648492575203252)
(7,0.6671791260162601)
(8,0.6922500589430894)
(9,0.728147551490515)
(10,0.7719864749322494)
(11,0.8175567890243902)
(12,0.8627245454200543)
(13,0.9024376513872019)
(14,0.9278605861585365)
(15,0.9370716494783198)
(16,0.9419975063482385)
(17,0.9437279045203253)
(18,0.9434962534688347)
(19,0.9420234507317073)
(20,0.9384741681233062)
(21,0.9309368872831308)
(22,0.919904599593496)
(23,0.8781241951219512)
            };
            
        \end{axis}
    \end{tikzpicture}
    \caption{Cost minimization G2V and V2G}
  \label{fig:HourlyLoadCostMin_v2gIL}
\end{subfigure}
\begin{subfigure}{0.33\textwidth}
 
\begin{tikzpicture}[scale=0.45, transform shape]
        \begin{axis}[
            xlabel = \Large Time (hour),
            title style={font=\large},
            xtick = {0,1,2,3,4,5,6,7,8,9,10,11,12,13,14,15,16,17,18,19,20,21,22,23},
            x tick label style={font=\large},
            ytick = {0, 0.2, 0.4, 0.6, 0.8, 1,1.2},
            y tick label style={font=\Large},
            ymin = 0,
            x = 0.42cm,
            legend pos = south east,
            xmajorgrids=false,
            grid style=dashed]
        \addplot[smooth, thick, color = blue, mark = *,style=dashed]
            plot coordinates {
                (0,0.81758257899729)
(1,0.7594078665311653)
(2,0.7131228883657181)
(3,0.6818801772357722)
(4,0.6601936653116531)
(5,0.6471568509485094)
(6,0.6484225630081301)
(7,0.6671175230352304)
(8,0.6921831266937669)
(9,0.7280820758807588)
(10,0.771910156504065)
(11,0.8174967323848239)
(12,0.8628398067005421)
(13,0.9025293252022283)
(14,0.9371358862127371)
(15,0.9699021759400142)
(16,0.9337487169101801)
(17,1.0028456747215753)
(18,1.007679794695122)
(19,0.9979741011924118)
(20,0.9743293061856368)
(21,0.9434269171414214)
(22,0.9207609478319784)
(23,0.8636606531165312)
                
            };
            \addplot[smooth, thick, color = purple, mark = *,style=solid]
            plot coordinates {
(0,0.04606345691056911)
(1,0.10420839024390242)
(2,1.1083178997289972e-05)
(3,0.029068434552845527)
(4,0.2034342337398374)
(5,0.21648560840108402)
(6,0.21522943157181573)
(7,0.19653110840108404)
(8,0.1714696795392954)
(9,0.13557386246612466)
(10,0.09175142005420053)
(11,0.046161615176151756)
(12,0.0010633694512195123)
(13,2.0731809688346882e-07)
(14,0.0016069098577235774)
(15,2.9221716734417347e-07)
(16,5.680255684281842e-08)
(17,1.5183669037940382e-07)
(18,0.0006321951422764227)
(19,0.0028574414905149052)
(20,0.002857467527100271)
(21,9.640871409214093e-08)
(22,0.0)
(23,0.0)
            };
            
\addplot[smooth, thick, color = orange, mark = *,style=dotted]
            plot coordinates {
               (0,0.0)
(1,0.0)
(2,0.0)
(3,0.0)
(4,0.0)
(5,0.0)
(6,0.0)
(7,0.0)
(8,0.0)
(9,0.0)
(10,0.0)
(11,0.0)
(12,2.731574132791328e-07)
(13,2.063587039295393e-06)
(14,2.8647551219512195e-06)
(15,4.589552560975609e-06)
(16,0.06222384681571815)
(17,3.434909701897019e-05)
(18,5.503700210027101e-06)
(19,6.953318428184282e-06)
(20,6.1546014701897014e-06)
(21,2.143975975609756e-06)
(22,6.38479508807588e-05)
(23,0.016243709620596205)
            };
        \end{axis}
    \end{tikzpicture}
    \caption{Intermediate case G2V and V2G}
  \label{fig:HourlyLoadCostMinInt_v2gIL}
\end{subfigure}
\begin{subfigure}{0.33\textwidth}
 
    \begin{tikzpicture}[scale=0.45, transform shape]
        \begin{axis}[
            xlabel = \Large Time (hour),
            title style={font=\large},
            xtick = {0,1,2,3,4,5,6,7,8,9,10,11,12,13,14,15,16,17,18,19,20,21,22,23},
            x tick label style={font=\large},
            ytick = {0, 0.2, 0.4, 0.6, 0.8, 1,1.2},
            y tick label style={font=\Large},
            ymin = 0,
            x = 0.42cm,
            legend pos = south east,
            legend style={at={(0.98,0.3)}},
            xmajorgrids=false,
            grid style=dashed]
        \addplot[smooth, thick, color = blue, mark = *,style=dashed]
            plot coordinates {
(0,0.8175891756097561)
(1,0.7585998752710027)
(2,0.7131220312968866)
(3,0.6814426283672087)
(4,0.6591154539295393)
(5,0.6459635819783198)
(6,0.6484979796747968)
(7,0.6654887935636857)
(8,0.6909667401084011)
(9,0.7281260697831977)
(10,0.7719952967479675)
(11,0.8204824593495935)
(12,0.8679179234417345)
(13,0.9027671748644986)
(14,0.9388650291327912)
(15,0.9702915483739837)
(16,0.9909655325474254)
(17,1.0028769718436608)
(18,1.0078467602371273)
(19,0.9981312311653115)
(20,0.9747917592818428)
(21,0.943428353320779)
(22,0.9208188921225938)
(23,0.878426656504065)
                
            };
            \addplot[smooth, thick, color = purple, mark = *,style=solid]
            plot coordinates {
                (0,0.046060104200542004)
(1,0.009679346273712737)
(2,3.2559588888888886e-09)
(3,0.0028255979200542004)
(4,0.09574088888888888)
(5,0.09883820799457993)
(6,0.21514943699186992)
(7,0.03413081890243902)
(8,0.043703207046070454)
(9,0.1355218096205962)
(10,0.09497961517615176)
(11,0.20294173780487806)
(12,0.2468228340108401)
(13,0.012926378658536585)
(14,0.08265642140921409)
(15,0.01615582086720867)
(16,0.00182791054200542)
(17,1.1826474051490516e-07)
(18,0.006178678109756098)
(19,0.008520628590785909)
(20,0.02300889315718157)
(21,3.3028355691056915e-07)
(22,2.31890593495935e-08)
(23,0.0)
                
            };
            
\addplot[smooth, thick, color = orange, mark = *,style=dotted]
            plot coordinates {
                (0,0.0)
(1,0.0)
(2,8.410470054200542e-08)
(3,0.0)
(4,0.0)
(5,0.0)
(6,0.0)
(7,0.0)
(8,0.0)
(9,0.0)
(10,0.0)
(11,0.0)
(12,0.0)
(13,0.0)
(14,0.0)
(15,0.0)
(16,1.0901226422764229e-07)
(17,2.3789160094850946e-09)
(18,0.0)
(19,0.0)
(20,0.0)
(21,0.0)
(22,3.514568299457994e-08)
(23,2.870595731707317e-06)
                
            };
            
        \addlegendentry{Electricity generation excluding EVs}
          \addlegendentry{Electricity generated for charging EVs}
        \addlegendentry{V2G power profile} 
        
        \end{axis}
    \end{tikzpicture}
    \caption{Emission minimization G2V and V2G}
  \label{fig:HourlyLoadEmisMin_v2gIL}
\end{subfigure}
    \caption{Hourly load variations for Illinois Summer (Intermediate cases are indicated in Fig. \ref{fig:IL-SCPareto_a}). }
    \label{fig:HourlyLoadIL}
\label{fig:ILWT}
\end{figure*}
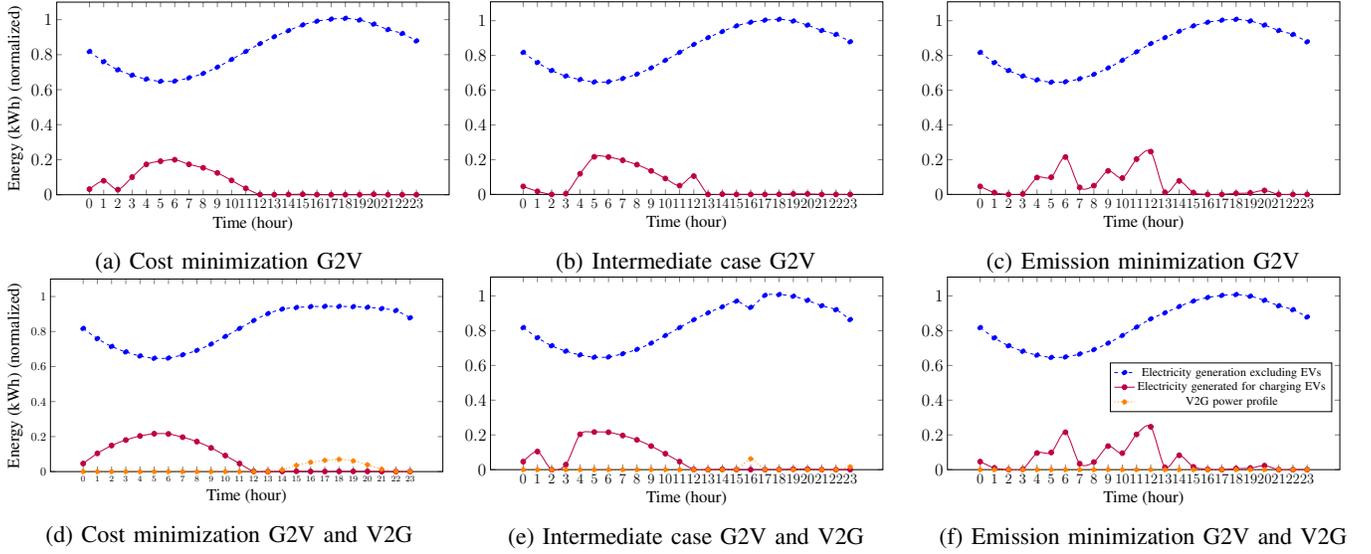

 \subsubsection{Hourly Load Variations}
 We now analyze the effects of hourly electricity consumption and marginal emission factors with respect to the trade-off between the cost and emission objectives. Under the cost minimization objective, we observe that the  peak of electricity generation for EV charging occurs when the generation excluding EVs is the lowest.
 For example, in Figs. \ref{fig:HourlyLoadCostMinIL} and \ref{fig:HourlyLoadCostMin}, when the conventional load is relatively low early in the morning, the electricity generated for charging EVs is high, while at the peak of the conventional load, the electricity generated for charging EVs is minimal. If the problem aims to minimize emission, the electricity is generated when the emission factor is relatively low (see Figs. \ref{fig:HourlyLoadEmisMinIL} and \ref{fig:HourlyLoadEmisMin}). This results in higher peaks of the grid; therefore, grid losses and, consequently, total generation cost increases. It is also important to consider the impacts of driving profiles and the initial battery state of charge. For example, in Fig. \ref{fig:HourlyLoadEmisMinIL}, although emission factors are higher in the early morning, EVs have to be charged in the early morning to satisfy their demand (see Fig. \ref{fig:drivingprofiles}). As expected, EV charging curves of intermediate cases are in between cost and emission minimization curves (see Figs. \ref{fig:HourlyLoadIntIL} and \ref{fig:HourlyLoadInt}).
 
Our findings can help decision-makers to design incentivization mechanisms for EV users depending on the preferred trade-off between cost and emission. For example, if cost-efficiency is the primary objective, then charging in the morning and feeding back to the grid in the evening can be promoted  {for peak shaving} (see Figs.  \ref{fig:HourlyLoadCostMin_v2gIL} and \ref{fig:HourlyLoadCostMin_v2g}).

\subsubsection{Effect of V2G}
We investigate how $\%$ changes in cost and emission are affected by the V2G concept. The results show that the bidirectional flow can decrease cost (e.g., Fig. \ref{fig:IL-SCPareto_c}). In the G2V-V2G case, EVs are charged during off-peak hours of electricity consumption and give energy back to the grid during peak hours to reduce the peak load of the power grid (see Fig. \ref{fig:HourlyLoadCostMin_v2g}). This flattens the load profile, resulting in a reduction of power losses, thereby decreasing the total generation cost. For the G2V-V2G case in Fig. \ref{fig:IL-SCPareto_a}, we even observe a cost reduction compared to the no-EV case through coordinated charging of EVs. According to Fig. \ref{fig:IL-SCPareto_d}, the integration of V2G does not significantly affect cost-emission trade-off since the demand profile is relatively straight in winter (see Fig. \ref{fig:SouthCarolinaConsumption}).

\subsubsection{Computational Performance}
To demonstrate the computational accuracy of our solution approach, we provide optimality gaps. For the 200-bus system in IL and the 2000-bus system in TX, we can provide approximately globally optimal solutions with an optimality gap below $0.1\%$ and  $0.6\%$, respectively. For the 500-bus system in SC, the optimality gaps are no more than 0.6\% in winter; however, in summer, the optimality gaps vary between $1.8 \%$ and $6.3 \%$. The authors in \cite{babaeinejadsarookolaee2019power} also note that the 500-bus case exhibits a significant SOCP optimality gap. Overall, our solution approach achieves relatively small percentage optimality gaps for all the instances we consider. We calculate computational times of obtaining Pareto frontiers of Illinois 200-bus, South Carolina 500-bus, and Texas 2000-bus cases as approximately  5 minutes, 21 minutes, and 42 minutes, respectively.

\section{Conclusions}\label{sec:MOPFConclusions}

This paper proposes a new mathematical formulation for the EV-aware MOPF Problem, and a corresponding solution approach based on SOCP. We consider a realistic OPF test case with three real datasets: hourly electricity consumption, hourly marginal emission factors, and EV driving profiles. We develop a convex programming framework based on SOCP to approximate globally optimal solutions. We conduct computational experiments to investigate the aggregate effect of the EV charging load on the grid. Computational experiments on different sizes of OPF instances from the PGLIB-OPF library show that our solution approach provides high quality solutions with small optimality gaps. 

We conduct our experiments at two different seasons and three different locations. Our results show that coordinated charging of EVs can decrease marginal $\text{CO}_2$ emissions significantly with either no change or minimal increase in cost. In addition, the integration of the V2G concept leads to cost savings. 

Our approach can have implications for real world management of EV charging and OPF. As an example, it is possible to  provide incentives to EV owners to charge their vehicles at times  that are cost and/or emission effective. Our methodology can provide insights to decision makers on when and how to incentivize EV charging.  

Further investigations can improve our results. In this study, we consider constant electricity prices.  
V2G might have greater cost reduction potential with variable electricity prices. In addition, we consider only one type of EV and all vehicles have the same constant charging rate. Including a mix of different electric vehicle models and charging rates would make a more realistic case, however this would complicate the optimization problem. Finally, our model requires hourly electricity consumption, EV charging demands, and the charging location of each EV to be known a priori. Considering the variability and stochasticity of these settings reveals itself as a challenging problem to be studied.

\bibliographystyle{ieeetr}
\bibliography{references}

\begin{thebibliography}{10}

\bibitem{AFDC}
AFDC, ``State alternative fuel and advanced vehicle laws and incentives: 2017
  year in review by alternative fuels data center.''
  \url{https://afdc.energy.gov/laws/}, 2017.

\bibitem{yuksel2016effect}
T.~Yuksel, M.-A.~M. Tamayao, C.~Hendrickson, I.~M. Azevedo, and J.~J. Michalek,
  ``Effect of regional grid mix, driving patterns and climate on the
  comparative carbon footprint of gasoline and plug-in electric vehicles in the
  united states,'' {\em Environmental Research Letters}, vol.~11, no.~4,
  p.~044007, 2016.

\bibitem{clement2009impact}
K.~Clement-Nyns, E.~Haesen, and J.~Driesen, ``The impact of charging plug-in
  hybrid electric vehicles on a residential distribution grid,'' {\em IEEE
  Trans. on Power Syst.}, vol.~25, no.~1, pp.~371--380, 2009.

\bibitem{gan2012optimal}
L.~Gan, U.~Topcu, and S.~H. Low, ``Optimal decentralized protocol for electric
  vehicle charging,'' {\em IEEE Trans. on Power Syst.}, vol.~28, no.~2,
  pp.~940--951, 2012.

\bibitem{masoum2012distribution}
M.~A. Masoum, P.~S. Moses, and S.~Hajforoosh, ``Distribution transformer stress
  in smart grid with coordinated charging of plug-in electric vehicles,'' in
  {\em 2012 IEEE PES Innovative Smart Grid Technologies (ISGT)}, pp.~1--8,
  IEEE, 2012.

\bibitem{judd2008evaluation}
S.~L. Judd and T.~J. Overbye, ``An evaluation of phev contributions to power
  system disturbances and economics,'' in {\em 2008 40th North American Power
  Symposium}, pp.~1--8, IEEE, 2008.

\bibitem{acha2010effects}
S.~Acha, T.~C. Green, and N.~Shah, ``Effects of optimised plug-in hybrid
  vehicle charging strategies on electric distribution network losses,'' in
  {\em IEEE PES T\&D 2010}, pp.~1--6, IEEE, 2010.

\bibitem{tang2016model}
W.~Tang and Y.~J.~A. Zhang, ``A model predictive control approach for
  low-complexity electric vehicle charging scheduling: Optimality and
  scalability,'' {\em IEEE Trans. on Power Syst.}, vol.~32, no.~2,
  pp.~1050--1063, 2016.

\bibitem{6847105}
N.~{Chen}, C.~W. {Tan}, and T.~Q.~S. {Quek}, ``Electric vehicle charging in
  smart grid: Optimality and valley-filling algorithms,'' {\em IEEE Journal of
  Selected Topics in Signal Processing}, vol.~8, no.~6, pp.~1073--1083, 2014.

\bibitem{fan2017admm}
H.~Fan, C.~Duan, C.-K. Zhang, L.~Jiang, C.~Mao, and D.~Wang, ``Admm-based
  multiperiod optimal power flow considering plug-in electric vehicles
  charging,'' {\em IEEE Trans. on Power Syst.}, vol.~33, no.~4, pp.~3886--3897,
  2017.

\bibitem{shi2018model}
Y.~Shi, H.~D. Tuan, A.~V. Savkin, T.~Q. Duong, and H.~V. Poor, ``Model
  predictive control for smart grids with multiple electric-vehicle charging
  stations,'' {\em IEEE Trans. on Smart Grid}, vol.~10, no.~2, pp.~2127--2136,
  2018.

\bibitem{7384482}
R.~{Azizipanah-Abarghooee}, V.~{Terzija}, F.~{Golestaneh}, and A.~{Roosta},
  ``Multiobjective dynamic optimal power flow considering fuzzy-based smart
  utilization of mobile electric vehicles,'' {\em IEEE Trans. on Industrial
  Informatics}, vol.~12, no.~2, pp.~503--514, 2016.

\bibitem{6485952}
N.~Chen, T.~Q. Quek, and C.~W. Tan, ``Optimal charging of electric vehicles in
  smart grid: Characterization and valley-filling algorithms,'' in {\em 2012
  IEEE Third International Conference on Smart Grid Communications
  (SmartGridComm)}, pp.~13--18, 2012.

\bibitem{yang2014improved}
J.~Yang, L.~He, and S.~Fu, ``An improved pso-based charging strategy of
  electric vehicles in electrical distribution grid,'' {\em Applied Energy},
  vol.~128, pp.~82--92, 2014.

\bibitem{valentine2011intelligent}
K.~Valentine, W.~G. Temple, and K.~M. Zhang, ``Intelligent electric vehicle
  charging: Rethinking the valley-fill,'' {\em Journal of Power Sources},
  vol.~196, no.~24, pp.~10717--10726, 2011.

\bibitem{alonso2014optimal}
M.~Alonso, H.~Amaris, J.~G. Germain, and J.~M. Galan, ``Optimal charging
  scheduling of electric vehicles in smart grids by heuristic algorithms,''
  {\em Energies}, vol.~7, no.~4, pp.~2449--2475, 2014.

\bibitem{yang2014self}
Z.~Yang, K.~Li, Q.~Niu, Y.~Xue, and A.~Foley, ``A self-learning tlbo based
  dynamic economic/environmental dispatch considering multiple plug-in electric
  vehicle loads,'' {\em Journal of Modern Power Systems and Clean Energy},
  vol.~2, no.~4, pp.~298--307, 2014.

\bibitem{yang2015computational}
Z.~Yang, K.~Li, and A.~Foley, ``Computational scheduling methods for
  integrating plug-in electric vehicles with power systems: A review,'' {\em
  Renewable and Sustainable Energy Reviews}, vol.~51, pp.~396--416, 2015.

\bibitem{zohrizadeh2020survey}
F.~Zohrizadeh, C.~Josz, M.~Jin, R.~Madani, J.~Lavaei, and S.~Sojoudi, ``A
  survey on conic relaxations of optimal power flow problem,'' {\em European
  journal of operational research}, vol.~287, no.~2, pp.~391--409, 2020.

\bibitem{li2012optimization}
N.~Li, L.~Gan, L.~Chen, and S.~H. Low, ``An optimization-based demand response
  in radial distribution networks,'' in {\em 2012 IEEE Globecom Workshops},
  pp.~1474--1479, IEEE, 2012.

\bibitem{gopalakrishnan2013global}
A.~Gopalakrishnan, A.~U. Raghunathan, D.~Nikovski, and L.~T. Biegler, ``Global
  optimization of multi-period optimal power flow,'' in {\em 2013 American
  Control Conference}, pp.~1157--1164, IEEE, 2013.

\bibitem{jabr2014minimum}
R.~A. Jabr, ``Minimum loss operation of distribution networks with photovoltaic
  generation,'' {\em IET Renewable Power Generation}, vol.~8, no.~1,
  pp.~33--44, 2014.

\bibitem{huang2016sufficient}
S.~Huang, Q.~Wu, J.~Wang, and H.~Zhao, ``A sufficient condition on convex
  relaxation of ac optimal power flow in distribution networks,'' {\em IEEE
  Trans. on Power Syst.}, vol.~32, no.~2, pp.~1359--1368, 2016.

\bibitem{kocuk2016strong}
B.~Kocuk, S.~S. Dey, and X.~A. Sun, ``Strong {SOCP} relaxations for the optimal
  power flow problem,'' {\em Oper. Res.}, vol.~64, no.~6, pp.~1177--1196, 2016.

\bibitem{birchfield2016grid}
A.~B. Birchfield, T.~Xu, K.~M. Gegner, K.~S. Shetye, and T.~J. Overbye, ``Grid
  structural characteristics as validation criteria for synthetic networks,''
  {\em IEEE Trans. on Power Syst.}, vol.~32, no.~4, pp.~3258--3265, 2016.

\bibitem{eia.gov2020}
EIA, ``Hourly electricy consumption data by united states energy information
  administration.''
  \url{http://https://www.eia.gov/beta/electricity/gridmonitor/}, 2020.

\bibitem{cedm2020}
I.~Azevedo, P.~Donti, N.~Horner, G.~Schivley, K.~Siler-Evans, and P.~Vaishnav,
  ``Electricity marginal factor estimates by center for climate and energy
  decision making pittsburgh: Carnegie mellon university.''
  \url{http://cedmcenter.org}, 2019.

\bibitem{NHTS}
NHTS, ``2017 national household travel survey by united states department of
  transportation, federal highway administration.''
  \url{https://nhts.ornl.gov}, 2017.

\bibitem{INSIDEEVs}
INSIDEEVs, ``Nissan leaf: "the world's most popular ev with 450,000 sold since
  2010.'' \url{https://insideevs.com/news/393890/nissan-leaf-sales-450000},
  2020.

\bibitem{nissanleaf}
NISSAN, ``Nissan leaf range \& batteries.''
  \url{https://www.nissanusa.com/vehicles/electric-cars/leaf/features/range-charging-battery.html},
  2020.

\bibitem{fueleco}
EPA, ``Nissan leaf's average energy consumption by united states environmental
  protection agency.'' \url{https://fueleconomy.gov/}, 2020.

\bibitem{richardson2011optimal}
P.~Richardson, D.~Flynn, and A.~Keane, ``Optimal charging of electric vehicles
  in low-voltage distribution systems,'' {\em IEEE Transactions on Power
  Systems}, vol.~27, no.~1, pp.~268--279, 2011.

\bibitem{afdcfuel}
AFDC, ``Average fuel economy by major vehicle category.''
  \url{https://afdc.energy.gov/data/10310}, 2020.

\bibitem{babaeinejadsarookolaee2019power}
S.~Babaeinejadsarookolaee, A.~Birchfield, R.~D. Christie, C.~Coffrin,
  C.~DeMarco, R.~Diao, M.~Ferris, S.~Fliscounakis, S.~Greene, R.~Huang, {\em
  et~al.}, ``The power grid library for benchmarking ac optimal power flow
  algorithms,'' {\em arXiv preprint arXiv:1908.02788}, 2019.

\end{thebibliography}

\end{document}